\documentclass[12pt]{amsart}

\usepackage{amsmath,amssymb}

\usepackage[all]{xy}
\CompileMatrices

\numberwithin{equation}{section}
\theoremstyle{plain}
\newtheorem{theorem}{Theorem}[section]
\newtheorem{proposition}[theorem]{Proposition}
\newtheorem{lemma}[theorem]{Lemma}
\newtheorem{corollary}[theorem]{Corollary}

\theoremstyle{definition}
\newtheorem{definition}[theorem]{Definition}
\newtheorem{examples}[theorem]{Examples} 
\newtheorem{example}[theorem]{Example}   

\newtheorem{remarks}[theorem]{Remarks}  
\newtheorem{remark}[theorem]{Remark} 

\newcommand{\bd}{\begin{definition}}
\newcommand{\ed}{\end{definition}}
\newcommand{\bt}{\begin{theorem}}
\newcommand{\et}{\end{theorem}}
\newcommand{\bp}{\begin{proposition}}
\newcommand{\ep}{\end{proposition}}
\newcommand{\bl}{\begin{lemma}}
\newcommand{\el}{\end{lemma}}
\newcommand{\bc}{\begin{corollary}}
\newcommand{\ec}{\end{corollary}}
\newcommand{\bpf}{\begin{proof}}
\newcommand{\epf}{\end{proof}}
\newcommand{\br}{\begin{remarks}\begin{flushleft}\end{flushleft}\nopagebreak} 
\newcommand{\er}{\end{remarks}}
\newcommand{\bex}{\begin{examples}\begin{flushleft}\end{flushleft}\nopagebreak} 
\newcommand{\eex}{\end{examples}}

\newcommand{\be}{\begin{enumerate}}
\newcommand{\ee}{\end{enumerate}}
\newcommand{\bi}{\begin{itemize}}
\newcommand{\ei}{\end{itemize}}

\newcommand{\sumsub}[1]{\sum_{\substack{\rule{0pt}{8pt}#1\rule{0pt}{8pt}}}} 

\newcommand{\Ualpha}                                                           
{\mbox{{$\alpha$}\hspace*{-8.3pt}\raisebox{-5pt}{$\scriptstyle{\smile}$}}}

\newcommand{\abs}[1]{\lvert#1\rvert}   

\newcommand{\rad}{\mathrm{rad}}
\newcommand{\tr}{\mathrm{tr}}
\newcommand{\ipart}[1]{\lfloor \frac{#1}{2} \rfloor} 
\newcommand{\ipartn}{\ipart{n}}  

\newcommand{\inv}[1]{\overline{#1}}

\newcommand{\Des}{\mathrm{Des}}
\newcommand{\Sh}{\mathrm{Sh}}

\newcommand{\Peak}{\mathrm{Peak}}

\newcommand{\Span}{\mathrm{Span}}
\renewcommand{\ker}{\mathrm{Ker}} 

\newcommand{\codim}{\mathrm{codim}\,}
\newcommand{\End}{\textrm{End}}

\newcommand{\Z}{\mathbb{Z}}

\newcommand{\I}{\mathfrak{I}}
\newcommand{\ppp}[1]{\mathfrak{P}_{#1}}
\newcommand{\pppn}{\ppp{n}}
\newcommand{\pppnj}{\pppn^j}

\newcommand{\TO}{\Bar{O}}

\newcommand{\pint}{p^{0}}

\newcommand{\hwp}{\widehat{\wp}}

\newcommand{\Sol}[1]{\mathit{Sol}(#1)}
\newcommand{\SolB}{\Sol{B_n}}
\newcommand{\SolA}{\Sol{A_{n-1}}}

\newcommand{\sol}[1]{\mathit{s}(#1)}
\newcommand{\sB}{\sol{B_n}}

\newcommand{\hsol}[1]{\widehat{\mathit{s}}(#1)}
\newcommand{\hsB}{\hsol{B_n}}

\newcommand{\field}{\Bbbk} 

\newcommand{\equal}[1]{
{\stackrel{{\textstyle #1}}{\ {\textstyle =}\ } }}
\newcommand{\onto}{\twoheadrightarrow}

\newcommand{\inc}{\hookrightarrow}

\begin{document}
\title[The peak algebra]{New results on the peak algebra} 
\author{Marcelo Aguiar}
\address{Department of Mathematics\\
         Texas A\&M University\\
         College Station\\
         TX  77843\\
         USA}
\email{maguiar@math.tamu.edu}
\urladdr{http://www.math.tamu.edu/$\sim$maguiar}

\author{Kathryn Nyman}
\address{Department of Mathematics\\
         Texas A\&M University\\
         College Station\\
         TX  77843\\
         USA}
\email{nyman@math.tamu.edu}
\urladdr{http://www.math.tamu.edu/$\sim$kathryn.nyman}

\author{Rosa Orellana}
\address{Department of Mathematics\\
        Dartmouth College \\
        Hanover, NH 03755\\
        USA}
\email{Rosa.C.Orellana@Dartmouth.EDU}
\urladdr{http://www.math.dartmouth.edu/$\sim$orellana/}

\thanks{Aguiar supported in part by NSF grant DMS-0302423}
\thanks{Orellana supported in part by the Wilson Foundation}

\date{June 9, 2004}

\begin{abstract}
The peak algebra $\pppn$ is a unital subalgebra of the symmetric group algebra, linearly spanned
by sums of permutations with a common set of peaks.  By exploiting the combinatorics of {\em sparse} subsets of $[n-1]$ (and of certain classes of compositions of $n$ called {\em almost-odd}  and {\em thin}), we construct
three new linear bases of $\pppn$. We  discuss two peak analogs of the first Eulerian idempotent and construct a basis of semi-idempotent elements
for the peak algebra. We use these bases to describe the Jacobson radical of $\pppn$ and to characterize
the elements of $\pppn$ in terms of the canonical action of  the symmetric groups on the
tensor algebra of a vector space. We define a chain of ideals $\pppnj$ of $\pppn$, $j=0,\ldots,\ipartn$,
such that $\pppn^{0}$ is the linear span of sums of permutations with a common set of interior peaks
and $\pppn^{\ipartn}$ is the peak algebra. We extend the above results to  $\pppnj$, generalizing
results of Schocker (the case $j=0$). \end{abstract}

\maketitle


\section*{Introduction}\label{S:intro}

A descent of a permutation $\sigma\in S_n$ is a position $i$ for which
$\sigma(i)>\sigma(i+1)$, while a peak is a position $i$ for which $\sigma(i-1)<\sigma(i)>\sigma(i+1)$. 

One aspect of the algebraic theory of peaks was
initiated by Stembridge~\cite{Ste}, another by Nyman~\cite{Nym}.
The peak algebra $\pppn$ was introduced in~\cite{ABN}. It is a unital subalgebra of the group algebra of the symmetric group $S_n$, obtained as the linear span of sums of permutations with a common set of peaks. The construction is analogous to that of the descent algebra of $S_n$, denoted $\SolA$, which is obtained as the linear span of sums of permutations with a common  set of descents. $\pppn$ is a subalgebra of $\SolA$.

The descent algebra has been the object of numerous works; for a recent survey see~\cite{Sch2}. The peak algebra, or closely related objects, has been studied in~\cite{ABN,BHT,Hsi,Sch}, from different perspectives.

 The descent algebra construction, due to Solomon, can be extended to all finite Coxeter groups~\cite{Sol}. Let $B_n$ be the group of signed permutations: $B_n=S_n\ltimes\Z_2^n$, and
\[\varphi:B_n\to S_n\]
 the canonical projection (the map that forgets the signs).
A basic observation of~\cite{ABN} is that this map sends the descent algebra of $B_n$, denoted $\SolB$, onto the peak algebra $\pppn$. This allows us to derive properties of the peak algebra from known properties of the descent algebra of $B_n$. This point of view is emphasized again in this work.

\subsection*{Notation} We write $[m,n]:=\{m,m+1,\ldots,n\}$ and $[n]:=[1,n]$. $\Z$ is the set of integers.
A subset $F$ of $\Z$ is {\em sparse} if it does not contain consecutive integers:
for any $i,j\in F$, $\abs{i-j}\neq 1$.
The number of sparse subsets on $[n-1]$ is the Fibonacci number $f_n$, defined by
\[f_0=f_1=1 \text{ and }f_n=f_{n-1}+f_{n-2} \text{ for }n\geq 2\,.\]
Unless otherwise stated, $F$, $G$, and $H$ denote sparse subsets of $[n-1]$.

For any $i\in\Z$ and  $J\subseteq\Z$, we let $J+i:=\{j+i\mid j\in J\}$.
We use mostly $i=\pm 1$.

Given a (signed or ordinary) permutation $\sigma$, we let $\sigma(0)=0$ and define
\[\Des(\sigma):=\{i\in[n-1] \mid \sigma(i)>\sigma(i+1)\},\quad
\Peak(\sigma):=\{i\in[n-1] \mid \sigma(i-1)<\sigma(i)>\sigma(i+1)\}\]
if $\sigma\in S_n$, and
\[\Des(\sigma):=\{i\in[0,n-1] \mid \sigma(i)>\sigma(i+1)\}\]
if $\sigma\in B_n$.  Note that a signed permutation may have a descent at $i=0$
(if $\sigma(1)<0$) and an ordinary permutation may have a peak at $i=1$ (if $\sigma(1)>\sigma(2)$). If  $\sigma\in S_n$, $\Des(\sigma)$ is a subset of $[n-1]$ and $\Peak(\sigma)$ is a sparse subset of $[n-1]$; if $\sigma\in B_n$, $\Des(\sigma)$ is a subset of $[0,n-1]$.  

We work over a field $\field$ of characteristic different from $2$. 

The descent algebra $\SolA$ is the subspace of $\field S_n$  linearly spanned by the elements
\[Y_I:=\sum_{\sigma\in S_n,\,\Des(\sigma)=I} \sigma\,,\]
or by the elements
\[X_I:=\sum_{\sigma\in S_n,\,\Des(\sigma)\subseteq I} \sigma\,,\]
as $I$ runs over the subsets of $[n-1]$.
The peak algebra $\pppn$ is the subspace of $\field S_n$  linearly spanned by the elements
\[P_F:=\sum_{\sigma\in S_n,\,\Peak(\sigma)=F} \sigma\,,\]
as $F$ runs over the sparse subsets of $[n-1]$. The descent algebra $\SolB$ is the subspace of $\field B_n$  linearly spanned by the elements
\[Y_J:=\sum_{\sigma\in B_n,\,\Des(\sigma)=J} \sigma\,,\]
or by the elements
\[X_J:=\sum_{\sigma\in B_n,\,\Des(\sigma)\subseteq J} \sigma\,,\]
as $J$ runs over the subsets of $[0,n-1]$. 

It is sometimes convenient to index  basis elements of $\SolA$ by compositions of $n$ and basis elements of $\SolB$ by {\em pseudocompositions} of $n$: integer sequences $(b_0,b_1,\ldots,b_k)$ such that $b_0\geq 0$, $b_i>0$, and $b_0+b_1+\cdots+b_k=n$ (see Section~\ref{S:compositions}).

In Section~\ref{S:basis}, $p_j$ denotes a certain element of the peak algebra,
but in Section~\ref{S:Lie} the same symbol is used for Lie polynomials.

\subsection*{Contents} In Section~\ref{S:bases}, we construct three bases ($Q$, $O$, and $\TO$) of the
peak algebra and describe how they relate to each other. Two different partial orders on the set of sparse subsets of $[n-1]$ play a crucial role here. Section~\ref{S:compositions} continues the study of the combinatorics of sparse subsets, by introducing two closely related classes of compositions (thin and almost-odd). One of the partial orders on sparse subsets corresponds to refinement of thin compositions, the other to refinement of almost-odd compositions (Lemmas~\ref{L:sparse-thin} and~\ref{L:sparse-almostodd}).
Basis elements of the peak algebra may be indexed by either sparse subsets, thin compositions, or almost-odd compositions; the most convenient choice depending on the situation.

A chain of ideals $\pppn^j$, $j=0,\ldots,\ipartn$, of the peak algebra is introduced in Section~\ref{S:chains}. 
The ideal at the bottom of the chain, $\pppn^0$, is the {\em peak ideal} of~\cite{ABN}. It is the linear span of sums of permutations with a common set of
interior peaks. This is the object studied in~\cite{BHT,Hsi,Nym,Sch}. Our results recover several known results for $\pppn^0$, and extend them to the ideals $\pppn^j$ and the peak algebra $\pppn$.
This chain of ideals is the image of a chain of ideals of $\SolB$ under the map $\varphi$ (Proposition~\ref{P:ideals}).

In Section~\ref{S:radical} we study the (Jacobson) radical of the peak algebra.
The radical  of the descent algebra of an arbitrary finite
Coxeter group was described by Solomon~\cite[Theorem 3]{Sol}; see also~\cite[Theorem 1.1]{GReu} for the case of type A and~\cite[Corollary 2.13]{B92} for the case of type B.  As $(a_1,\ldots,a_k)$ runs over all compositions of $n$ and $s$ over all permutations of $[k]$,  the elements
\[X_{(a_1,\ldots,a_k)}-X_{(a_{s(1)},\ldots,a_{s(k)})}\]
 linearly span $\rad(\SolA)$, while
$\rad(\SolB)$ is linearly spanned by the elements
\[X_{(b_0,b_1,\ldots,b_k)}-X_{(b_0,b_{s(1)},\ldots,b_{s(k)})}\]
as $(b_0,b_1,\ldots,b_k)$ runs over all pseudocompositions of $n$ and $s$ over all permutations of $[k]$.
In Theorem~\ref{T:radical} we obtain a similar result for the radical of $\pppn$: 
$\rad(\pppn)$ is linearly spanned by the elements
\[Q_{(b_0,b_1,\ldots,b_k)}-Q_{(b_0,b_{s(1)},\ldots,b_{s(k)})}\]
as $(b_0,b_1,\ldots,b_k)$ runs over all almost-odd compositions of $n$ and $s$ over all permutations of $[k]$ (a similar result holds for the bases $O$ and $\TO$ as well). It follows that the codimension of the radical is the number of
almost-odd partitions of $n$ (Corollary~\ref{C:radical}). We also obtain similar
descriptions for the intersection of the radical with the ideals $\pppn^j$. 
The case $j=0$ recovers a result of Schocker on the radical of the peak ideal~\cite[Corollary 10.3]{Sch}.

Section~\ref{S:convolution} discusses the {\em external} structure on the direct sum of the peak algebras. This is a product on the space $\ppp{}=\oplus_{n\geq 0}\pppn$ which corresponds to the convolution  product of endomorphisms of the tensor
algebra $T(V)=\oplus_{n\geq 0}V^{\otimes n}$ via the canonical action of $S_n$ on $V^{\otimes n}$. The connection with the convolution product on $\Sol{B}=\oplus_{n\geq 0}\SolB$ is explained, and then used to derive properties of the convolution product on $\ppp{}$ from
properties of the convolution product on $\Sol{B}$, which is simpler to analyze.
Proposition~\ref{P:peakfree} states that the bases $Q$, $O$, and $\TO$ are
multiplicative with respect to the convolution product. It follows that $\ppp{}^0=\oplus_{n\geq 0}\pppn^0$ is a free algebra (with respect to the convolution product) with one generator for each odd degree (a result known from~\cite{BMSW,Hsi,Sch}) and that 
$\ppp{}$ is  free as a right module over $\ppp{}^0$, with one generator for each even degree.

Let $L(V)$ be the free Lie algebra generated by $V$. It is the subspace of primitive elements of the tensor algebra $T(V)$. The elements of $L(V)$ are called Lie polynomials and products of these are called Lie monomials. 
 The {\em first Eulerian idempotent} is a certain element of $\SolA$ which projects  the homogeneous component of degree $n$ of $T(V)$ onto the homogeneous component of degree $n$ of $L(V)$, 
via the canonical action of the symmetric groups  on the tensor algebra.
 The Eulerian idempotents have been thoroughly studied~\cite[Section 4.5]{Lod}, \cite[Chapter 3]{Reu}. In Section~\ref{S:basis} we discuss two {\em peak analogs} of the first Eulerian idempotent, $\rho_{(n)}$ and $\rho_{(0,n)}$. The latter was introduced by Schocker~\cite[Section 7]{Sch}. The former is idempotent when $n$ is even, the latter when $n$ is odd. 
 We describe these elements explicitly in terms of sums of permutations with a common number of peaks and show that they are images under $\varphi$ of elements introduced by Bergeron and Bergeron (Theorem~\ref{T:etorho}).
 We use them as the building blocks for a multiplicative basis of $\pppn$ consisting of semiidempotents elements (Corollary~\ref{C:rhobasis}).
 The idempotents $\rho_{(0,n)}$ ($n$ odd) project onto the  odd components of $L(V)$, while the idempotents $\rho_{(n)}$ ($n$ even) project onto the subalgebra of $T(V)$ generated by the even components of $L(V)$ (Lemma~\ref{L:rho-lie}). The elements $\rho_{(n)}$ and  $\rho_{(0,n)}$ belong
 to a commutative semisimple subalgebra of $\pppn$ introduced in~\cite[Section 6]{ABN}. More information about this subalgebra is provided in Section~\ref{S:commutative}.

Section~\ref{S:Lie} contains our main results. The proofs rely on most of the preceding constructions. A classical result (Schur-Weyl duality) states that  if $\dim V\geq n$ then $\field S_n$ may be recovered as those endomorphisms of $V^{\otimes n}$ which
commute with the diagonal action of $GL(V)$.
Similarly, an important result of Garsia and Reutenauer characterizes
which elements of the group algebra $\field S_n$ belong to the descent algebra $\SolA$ in terms of their action on  Lie monomials
~\cite[Theorem 4.5]{GReu}: an element $\phi\in\field S_n$ belongs to $\SolA$ if and only if its action on an arbitrary Lie monomial $m$ yields a linear combination of Lie monomials each of which consists of a permutation of the factors of $m$; see~\eqref{E:action-GR}.  Schocker obtained a 
characterization for the elements of the peak ideal $\pppn^0$
in terms of the action on Lie monomials~\cite[Main Theorem 8]{Sch}: an element $\phi\in\SolA$ belongs to $\pppn^0$ if and only if its action annihilates any Lie monomial whose first factor is of even degree; see~\eqref{E:action-Schocker}.
We present a characterization for the elements of the
peak algebra $\pppn$ that is analogous to that of Garsia and Reutenauer, both in content and proof (Theorem~\ref{T:charpeak}).  Our result states that 
an element $\phi\in\field S_n$ belongs to $\pppn$ if and only if its action on an arbitrary Lie monomial $m$ in which all factors of even degree precede all factors of odd degree yields a linear combination of Lie monomials each of which consists of the even factors of $m$ (in the same order) followed by a permutation of the odd factors of $m$; see~\eqref{E:subspace}.
Furthermore, we provide
a characterization for the elements of
each ideal $\pppn^j$ that interpolates between Schocker's
characterization of the peak ideal and our characterization
of the peak algebra (Theorem~\ref{T:charpeakideal}). The action of  an element 
of $\pppn^j$ must in addition annihilate any Lie monomial $m$ as above in which the degree of the even part is larger than $2j$; see~\eqref{E:action-us}.

\subsection*{Acknowledgements} We thank Nantel Bergeron, Steve Chase, Sam Hsiao, and Swapneel Mahajan for useful comments.

\section{Bases of the peak algebra} \label{S:bases}

For any  subset $M\subseteq[n-1]$, let
\[\Bar{M}:=\{i\in[n-1] \mid \text{ either $i$ is in $M$ or both $i-1$ and $i+1$ are in $M$}\}\,.\]
In other words,
\[\Bar{M}=M\cup\bigl((M-1)\cap(M+1)\bigr)\,.\]
Note that
\begin{equation}\label{E:closure}
\Bar{\Bar{M}}=\Bar{M} \text{ \ and \ }(M\subseteq N\Rightarrow \Bar{M}\subseteq\Bar{N})\,.
\end{equation}
\bd\label{D:bases} For any sparse subset $F\subseteq[n-1]$, let
\begin{align}\label{E:QinP}
Q_F:= & \sum_{F\subseteq G} P_G\,,\\
\label{E:OinP}
O_F:= & \sum_{G\subseteq [n-1]\setminus F} P_G\,,\\
\label{E:TOinP}
\TO_F:= & \sum_{G\subseteq [n-1]\setminus \Bar{F}} P_G\,;
\end{align}
\ed
 in each case the sum being over sparse subsets $G$ of $[n-1]$. 
 For example, when $n=6$, 
\begin{align*}
Q_{\{1,3\}}  &=P_{\{1,3\}}+P_{\{1,3,5\}}\,,\\
O_{\{1,3\}}  &= P_{\emptyset} + P_{\{2\}} + P_{\{4\}} + P_{\{5\}}+ P_{\{2,4\}}+ P_{\{2,5\}}\,,\\
\TO_{\{1,3\}}  &=P_{\emptyset}  + P_{\{4\}} + P_{\{5\}} \,.
\end{align*} 

 View the collection of sparse subsets of $[n-1]$ as a poset under inclusion. All subsets of a sparse subset are again sparse; therefore, each interval of this poset is Boolean. Hence,~\eqref{E:QinP}
 is equivalent to
 \begin{equation}\label{E:PinQ}
P_F:=\sum_{F\subseteq G} (-1)^{\#G\setminus F}Q_G\,.
\end{equation}
Thus, as $F$ runs over the sparse subsets of $[n-1]$, the elements $Q_F$ form a linear basis of $\pppn$. The matrices relating the elements $P_G$ to the elements $O_F$ and $\TO_F$ are not triangular. However,  these elements also form linear bases of $\pppn$. This will be shown shortly (Corollary~\ref{C:OandTObasis}).

\begin{lemma}\label{L:Q-P} For any subset $M\subseteq [n-1]$,
\begin{equation}\label{E:Q-P}
\sumsub{G \text{ sparse}\\ G\subseteq [n-1]\setminus M} (-1)^{\#G}Q_G=
\sumsub{H \text{ sparse}\\ H\subseteq M} P_H\,.
\end{equation}
\end{lemma}
\begin{proof} We have
\[\sumsub{G \text{ sparse}\\ G\subseteq [n-1]\setminus M} \!\!\!\!(-1)^{\#G}Q_G
\equal{\eqref{E:QinP}}
\!\!\!\!\sumsub{G \text{ sparse}\\ G\subseteq [n-1]\setminus M} 
\sumsub{H \text{ sparse}\\G \subseteq H}(-1)^{\#G}P_H=
\sum_{H \text{ sparse}}
\Bigl(\sumsub{G\subseteq ([n-1]\setminus M)\cap H} (-1)^{\#G}\Bigr)P_H\,.\]
The inner sum is $1$ if $([n-1]\setminus M)\cap H=\emptyset$ and $0$ otherwise; ~\eqref{E:Q-P} follows.
\end{proof}
\begin{proposition}\label{P:OinQ} For any sparse subset $F\subseteq[n-1]$,
\begin{equation}\label{E:OinQ}
O_F=\sum_{G\subseteq F} (-1)^{\#G}Q_G\,.
\end{equation}
\end{proposition}
\begin{proof} Apply Lemma~\ref{L:Q-P} with $M=[n-1]\setminus F$.
\end{proof}

For each subset $J$ of $[0,n-1]$, let
$X_J=\sum_{\Des(\sigma)\subseteq J}\sigma$. As mentioned in the introduction, these elements form a basis
of  $\SolB$.

Let $\varphi:B_n\to S_n$ be the canonical map. In~\cite[Proposition 3.3]{ABN}, we showed that
for any $J\subseteq[0,n-1]$,
\begin{equation}\label{E:XinP}
\varphi(X_J)=2^{\# J}\cdot\!\!\!\!\sumsub{H \text{ sparse}\\ H\subseteq J\cup(J+1)} P_H\,.
\end{equation}

 \begin{proposition}\label{P:XinQ}
For any $J\subseteq[0,n-1]$,
\begin{equation}\label{E:XinQ}
\varphi(X_J)=2^{\# J}\cdot\!\!\!\!\sumsub{G \text{ sparse}\\ G\subseteq [n-1]\setminus\bigl(J\cup(J+1)\bigr)} (-1)^{\#G}Q_G\,.
\end{equation}
\end{proposition}
\begin{proof} Apply Lemma~\ref{L:Q-P} with $M=\bigl(J\cup(J+1)\bigr)\cap [n-1]$.
\end{proof}

Given sparse subsets $F$ and $G$ of $[n-1]$, define
\[F\preceq G \iff \Bar{F}\supseteq G\,.\]

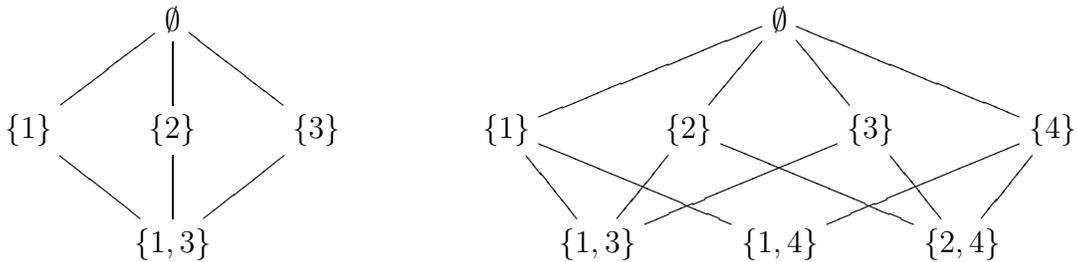
\begin{figure}[!h]
\[\xymatrix{
& \emptyset\ar@{-}[ld]\ar@{-}[d]\ar@{-}[rd] &\\
\{1\}\ar@{-}[rd] & \{2\}\ar@{-}[d] & \{3\}\ar@{-}[ld] \\
& \{1,3\} &} \qquad\qquad
\xymatrix@C-20pt{
& & & \emptyset\ar@{-}[llld]\ar@{-}[ld]\ar@{-}[rd]\ar@{-}[rrrd] & & & \\
\{1\}\ar@{-}[rd] \ar@{-}[rrrd] &  & \{2\}\ar@{-}[ld]  \ar@{-}[rrrd] &  & \{3\}\ar@{-}[llld] \ar@{-}[rd] & & \{4\}\ar@{-}[ld] \ar@{-}[llld] \\
& \{1,3\} & & \{1,4\} & & \{2,4\} &}
\]
\caption{Sparse subsets under $\preceq$}\label{F:sparse}
\end{figure}

\begin{lemma}\label{L:poset}
The  relation $\preceq$ is a partial order on the collection of sparse subsets of $[n-1]$.
\end{lemma}
\begin{proof} Suppose $F\preceq G$ and $G\preceq F$. Let $f=\max F$. Suppose $f\notin G$.
Then $f-1$ and $f+1\in G$, since $F\subseteq\Bar{G}$. Since $F$ is sparse, $f+1\notin F$. But then $f$ and $f+2\in F$, since $G\subseteq\Bar{F}$.
This contradicts the choice of $F$. Thus $f\in G$. Proceeding by induction, $F=G$.
This proves antisymmetry. Transitivity follows from~\eqref{E:closure}.
\end{proof}
The previous result may also be deduced from Lemma~\eqref{L:sparse-almostodd} .
The Hasse diagram of the poset of sparse subsets of $[n-1]$ under $\preceq$ are
shown in Figure~\ref{F:sparse}, for $n=4,5$.

\begin{proposition}\label{P:TOinQ} For any sparse subset $F\subseteq[n-1]$, 
\begin{equation}\label{E:TOinQ}
\TO_F=\sum_{F\preceq G} (-1)^{\#G}Q_G\,.
\end{equation}
\end{proposition}
\begin{proof} Let $J:=[0,n-1]\setminus \bigl(F\cup(F-1)\bigr)$. Then $J+1=[1,n]\setminus \bigl((F+1)\cup F\bigr)$.

On the other hand, 
\[\Bar{F}=F\cup\bigl((F-1)\cap(F+1)\bigr)=\bigl(F\cup(F-1)\bigr)\cap\bigl((F+1)\cup F\bigr)\,.\]
Therefore,
\[\bigl(J\cup(J+1)\bigr)\cap[n-1]=[n-1]\setminus \Bar{F}\,.\]
Combining~\eqref{E:TOinP} and~\eqref{E:XinP}  we deduce
\begin{equation*} \label{E:XtoTO}
\varphi(X_J)=2^{\# J}\cdot \TO_F\,.
\end{equation*}
Together with~\eqref{E:XinQ} this implies
\[\TO_F=\sumsub{G\text{ sparse}\\ G\subseteq [n-1]\setminus\bigl(J\cup(J+1)\bigr)} (-1)^{\#G}Q_G\,.\]
This sets~\eqref{E:TOinQ}, since by the above, $G\subseteq [n-1]\setminus\bigl(J\cup(J+1)\bigr)\iff G\subseteq\Bar{F} \iff F\preceq G$.
\end{proof}

\begin{corollary}\label{C:OandTObasis}
As $F$ runs over the sparse subsets of $[n-1]$, the elements $O_F$  form a linear basis of $\pppn$, and so do the elements $\TO_F$.
\end{corollary}
\begin{proof} Applying M\"obius inversion to~\eqref{E:OinQ} we obtain
\[Q_F=\sum_{G\subseteq F} (-1)^{\#G}O_G\,.\]
Let $\mu$ denote the M\"obius function of the poset of sparse subsets of $[n-1]$ under $\preceq$. Applying M\"obius inversion to~\eqref{E:TOinQ} we obtain
\[(-1)^{\#F}Q_F=\sum_{F\preceq G} \mu(F,G)\TO_G\,.\]
Since the elements $Q_F$ form a linear basis of $\pppn$, the same is true of
the elements $O_F$ and $\TO_F$.
\end{proof}
 The values $\mu(F,G)$ are products of Catalan numbers, see Remark~\ref{R:Catalan}.  
 Note that $\{P_F\}$, $\{Q_F\}$, $\{O_F\}$, and $\{\TO_F\}$ are integral bases of the peak algebra.

\section{Sparse subsets and  compositions} \label{S:compositions}

Let $n$ be a non-negative integer. An {\em ordinary composition} of $n$ is a sequence $\alpha=(a_1,\ldots,a_k)$
 of positive integers  such that  $a_1+\cdots +a_k=n$. 
  A {\em thin composition} of $n$ is an ordinary composition $\alpha$ of $n$  in which each $a_i$ is either $1$ or $2$.

A {\em pseudocomposition} of $n$ is a sequence $\beta=(b_0,b_1,\ldots,b_k)$ of integers such
that $b_0\geq 0$, $b_i\geq 1$ for $i\geq 1$, and $b_0+b_1+\cdots +b_k=n$. 
An {\em almost-odd composition} of $n$ is a pseudocomposition $\beta$ of $n$  in which
 $b_0\geq 0$ is even and $b_i\geq 1$ is odd for all $i\geq 1$.
 
We do not regard ordinary compositions as particular pseudocompositions. In particular, for  ordinary or thin compositions  $\alpha=(a_1,\ldots,a_k)$  we define the {\em number of parts} of $\alpha$ as
\[k(\alpha)=k\,,\]
 but for pseudo or almost-odd compositions $\beta=(b_0,b_1,\ldots,b_k)$  we define 
\begin{equation}\label{E:parts}
k(\beta)=k
\end{equation}
 (instead of $k+1$).

Pseudocompositions of $n$ are in bijection with subsets of $[0,n-1]$ via
\begin{equation}\label{E:pseudocomp-subset}
\beta=(b_0,b_1,\ldots,b_k)\mapsto J(\beta) := \{b_0, b_0+b_1, \ldots, b_0+ b_1+ \cdots
+b_{k-1}\}\,.
\end{equation}
Similarly,  compositions of $n$ are in bijection with subsets of $[n-1]$ via
\begin{equation*}\label{E:comp-subset}
\alpha=(a_1,\ldots,a_k)\mapsto I(\alpha) := \{a_1, a_1+a_2,\ldots, a_1+ a_2+ \cdots
+a_{k-1}\}\,.
\end{equation*}
 Under these bijections,
inclusion of subsets corresponds to refinement of compositions: $\beta'$ refines $\beta$ if and only if
$J(\beta)\subseteq J(\beta')$. We write $\beta\leq\beta'$ in this case. Note that
\[\#J(\beta)=k(\beta) \text{ \ and \ } \#I(\alpha)=k(\alpha)-1\,.\]
We use these correspondences to  label basis elements of $\SolB$ by pseudocompositions instead of subsets: given a pseudocomposition $\beta$ of $n$ we let $X_{\beta}:=X_{J(\beta)}$. Similarly, we may label basis elements of $\SolA$ by ordinary compositions of $n$.

There is a simple bijection between thin compositions of $n$ and sparse subsets of $[n-1]$.

\begin{lemma}\label{L:sparse-thin}Given a sparse subset $F$ of $[n-1]$, let $\tau_F$ be the unique ordinary composition of $n$
such that
\[I(\tau_F)=[n-1]\setminus F\,.\]
\begin{itemize}
\item[(i)] The composition $\tau_F$ is thin and
\[\#F=n-k(\tau_F)\,.\]
\item[(ii)] $F\mapsto \tau_F$ is a bijection between  sparse subsets of $[n-1]$ and thin compositions of $n$.
\item[(iii)]  Let $G$ be a sparse subset of $[n-1]$,  $\alpha$  an ordinary composition of $n$, and $I=I(\alpha)$. Then
\[G\subseteq [n-1]\setminus I \iff \alpha\leq\tau_G\,.\]
\item[(iv)]  For any sparse subsets $F$ and $G$ of $[n-1]$,
\[G \subseteq F \iff \tau_F \leq\tau_G\,.\]
\end{itemize}  
\end{lemma}
\begin{proof} Straightforward.
\end{proof}

\begin{figure}[!h]
\[\xymatrix@C-20pt{
& (1,1,1,1)\ar@{-}[ld]\ar@{-}[d]\ar@{-}[rd] &\\
(2,1,1)\ar@{-}[rd] & (1,2,1) & (1,1,2)\ar@{-}[ld] \\
& (2,2) &} \qquad
\xymatrix@C-40pt{
& & & (1,1,1,1,1)\ar@{-}[llld]\ar@{-}[ld]\ar@{-}[rd]\ar@{-}[rrrd] & & & \\
(2,1,1,1)\ar@{-}[rd] \ar@{-}[rrrd] &  & (1,2,1,1)  \ar@{-}[rrrd] &  & (1,1,2,1)\ar@{-}[llld]  & & (1,1,1,2)\ar@{-}[ld] \ar@{-}[llld] \\
& (2,2,1) & & (2,1,2) & & (1,2,2) &}
\]
\caption{Thin compositions under refinement}\label{F:thin}
\end{figure}
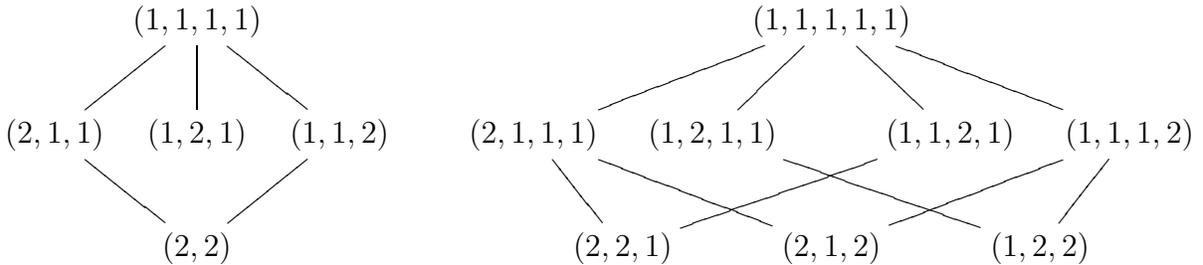

According to the lemma, the poset of sparse subsets of $[n-1]$ under reverse inclusion is isomorphic to
the poset of thin compositions of $n$ under refinement. The Hasse diagrams of the
latter are shown in Figure~\ref{F:thin}, for $n=4,5$. Comparison with Figure~\ref{F:sparse}
illustrates the correspondence of Lemma~\ref{L:sparse-thin} .

\medskip

There is also a bijection between almost-odd compositions of $n$ and sparse subsets of $[n-1]$.

\begin{lemma}\label{L:sparse-almostodd}Given a sparse subset $F$ of $[n-1]$, let $\gamma_F$ be the unique pseudocomposition of $n$
such that
\[J(\gamma_F)=[0,n-1]\setminus \bigl(F\cup(F-1)\bigr)\,.\]
\begin{itemize}
\item[(i)] The pseudocomposition $\gamma_F$ is almost-odd and
\[\#F=\frac{n-k(\gamma_F)}{2}\,.\]
\item[(ii)] $F\mapsto \gamma_F$ is a bijection between  sparse subsets of $[n-1]$ and almost-odd compositions of $n$.
\item[(iii)]  Let $G$ be a sparse subset of $[n-1]$,  $\beta$  a pseudocomposition of $n$, and $J=J(\beta)$. Then
\[G\subseteq [n-1]\setminus\bigl(J\cup(J+1)\bigr) \iff \beta\leq\gamma_G\,.\]
\item[(iv)]  For any sparse subsets $F$ and $G$ of $[n-1]$,
\[F\preceq G \iff  G\cup(G-1)\subseteq F\cup(F-1) \iff \gamma_F \leq\gamma_G\,.\]
\end{itemize}  
\end{lemma}
\begin{proof} We show (i). Since $F$ is sparse, it is a disjoint union of maximal subsets of the form $\{a,a+2,\ldots,a+2k\}$. It follows that 
 $F\cup(F-1)$ is a disjoint union of maximal intervals of the form
 $\{a-1,a,\ldots,a+2k-1,a+2k\}$. The difference between two consecutive elements of $J(\gamma_F)=[0,n-1]\setminus \bigl(F\cup(F-1)\bigr)$ is therefore odd (equal to $a+2k+1-a-2$). Consider the first element $a_0$ of $F$ and the corresponding interval $\{a_0-1,a_0,\ldots,a_0+2k_0-1,a_0+2k_0\}$. If $a_0=1$ then the first element of $J(\gamma_F)$ is $a_0+2k_0+1$ which is even. If $a_0=0$ then the first element of $J(\gamma_F)$ is $0$. This proves that $\gamma_F$ is almost-odd.
Also, $k(\gamma_F)=\#J(\gamma_F)=n-2\#F$, since $F$ and $F-1$ are disjoint and equinumerous.

Given an almost-odd composition $\gamma$, write $[0,n-1]\setminus J(\gamma)$ as a disjoint union of maximal intervals and delete every other element, starting with the first element of each interval. The result is a sparse subset of $[n-1]$.
This defines the inverse correspondence to $F\mapsto\gamma_F$, which proves (ii).

We show (iii). Refinement of pseudocompositions corresponds to inclusion of subsets via $J$. Therefore,
\begin{align*}
\beta\leq\gamma_G & \iff J(\beta)\subseteq J(\gamma_G) \iff G\cup(G-1)\subseteq [0,n-1]\setminus J\\
&\iff G\subseteq ([0,n-1]\setminus J)\cap\bigl([1,n]\setminus (J+1)\bigr)\\
&\iff G\subseteq [n-1]\setminus \bigl(J\cup(J+1)\bigr)\,.
\end{align*}

We show (iv). Let $\beta=\gamma_F$. Then $J=[0,n-1]\setminus\bigl(F\cup(F-1)\bigr)$.
The proof of (iii) shows that $\gamma_F \leq\gamma_G \iff G\cup(G-1)\subseteq F\cup(F-1)$.
The proof of Proposition~\ref{P:TOinQ} shows that $\bigl(J\cup(J+1)\bigr)\cap[n-1]=[n-1]\setminus \Bar{F}$.
Together with (iii) this says 
\[\gamma_F\leq \gamma_G \iff G\subseteq\Bar{F}\iff F\preceq G\,.\]
\end{proof}

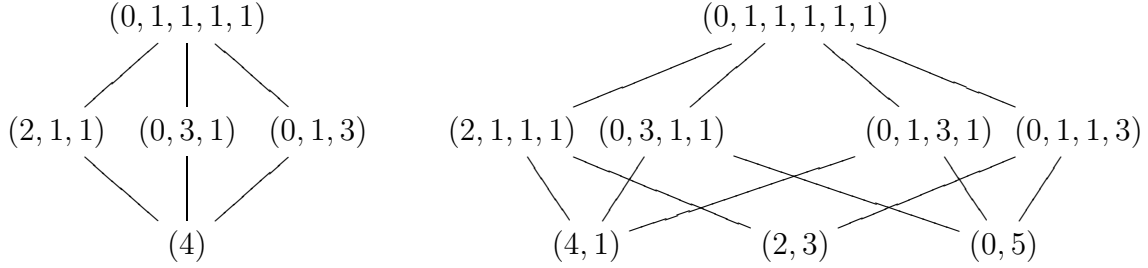
\begin{figure}[!h]
\[\xymatrix@C-30pt{
& (0,1,1,1,1)\ar@{-}[ld]\ar@{-}[d]\ar@{-}[rd] &\\
(2,1,1)\ar@{-}[rd] & (0,3,1)\ar@{-}[d] & (0,1,3)\ar@{-}[ld] \\
& (4) &} \qquad
\xymatrix@C-40pt{
& & & (0,1,1,1,1,1)\ar@{-}[llld]\ar@{-}[ld]\ar@{-}[rd]\ar@{-}[rrrd] & & & \\
(2,1,1,1)\ar@{-}[rd] \ar@{-}[rrrd] &  & (0,3,1,1)\ar@{-}[ld]  \ar@{-}[rrrd] &  & (0,1,3,1)\ar@{-}[llld] \ar@{-}[rd] & & (0,1,1,3)\ar@{-}[ld] \ar@{-}[llld] \\
& (4,1) & & (2,3) & & (0,5) &}
\]
\caption{Almost-odd compositions under refinement}\label{F:almostodd}
\end{figure}

According to the lemma, the poset of sparse subsets of $[n-1]$ under $\preceq$ is isomorphic to
the poset of almost-odd compositions of $n$ under refinement. The Hasse diagrams of the
latter are shown in Figure~\ref{F:almostodd}, for $n=4,5$. Comparison with Figure~\ref{F:sparse}
illustrates the correspondence of Lemma~\ref{L:sparse-almostodd} .

\begin{remark}\label{R:Catalan} 
The poset of almost-odd compositions of $n$ is isomorphic to the poset of
odd compositions of $n+1$ (add $1$ to the first part). 
It follows from~\cite[Exercise 52, Chapter 3]{St86} that the values of the M\"obius function of this poset are products of Catalan numbers. (The poset studied in this reference is the poset of odd compositions of $2m+1$. The poset of
odd compositions of $2m$ is a convex subset of the poset of odd compositions of $2m+1$: add a new part equal to $1$ at the end.) We thank Sam Hsiao for this reference.
\end{remark}

Combining the correspondences of Lemmas~\ref{L:sparse-thin} and~\ref{L:sparse-almostodd} results in a bijection between 
 thin compositions of $n$ and almost-odd compositions of $n$ 
 that we now describe.
 
 \begin{lemma}\label{L:thin-almostodd} Given an almost-odd composition $\gamma=(b_0,b_1,b_2,\ldots,b_k)$, let $\tau_\gamma$ be the  thin composition of $n$ given by
\[\tau_\gamma:=(\underbrace{2,\ldots,2}_{\frac{b_0}{2}},1,\underbrace{2,\ldots,2}_{\frac{b_1-1}{2}},1,\underbrace{2,\ldots,2}_{\frac{b_2-1}{2}},\ldots,
1,\underbrace{2,\ldots,2}_{\frac{b_{k}-1}{2}})\,.\]
\begin{itemize}
\item[(i)] $\gamma\mapsto \tau_\gamma$ is a bijection between   almost-odd compositions of $n$ and thin compositions of $n$ such that
\[k(\tau_\gamma)=\frac{n+k(\gamma)}{2}\,.\]
\item[(ii)] For any almost-odd compositions $\gamma$ and $\delta$ of $n$, we have
that $\tau_\gamma\leq \tau_\delta$ (refinement of thin compositions) if and only if
$\gamma\leq\delta$ (refinement of almost-odd compositions) and in addition
 $\delta$ is obtained by replacing each part $c$ of $\gamma$ by a sequence
 of parts $c_0,c_1,\ldots,c_i$ such that $c_0+c_1+\ldots+c_i=c$, $c_0\equiv c \mod 2$, $c_1=\ldots=c_{i-1}=1$, and $c_i$ is odd. (In particular, $i$ must be even.)
 \item[(iii)] The bijections $F\mapsto\tau_F$ and $F\mapsto\gamma_F$ of Lemmas~\ref{L:sparse-thin} and~\ref{L:sparse-almostodd} combine to give
 the bijection of (i), in the sense that $\tau_{\gamma_F}=\tau_F$.
\end{itemize}  
\end{lemma}
\begin{proof} Left to the reader.
\end{proof}
 
 For example, let $\gamma=(4,1,1)$ and $\delta=(0,3,1,1,1)$. Then
 \[\tau_\gamma=(2,2,1,1) \text{ \ and \ }\tau_\delta=(1,2,1,1,1)\,.\] Note that
 $\delta$ refines $\gamma$ but $\tau_\delta$ does not refine $\tau_\gamma$.
 In passing from $\gamma$ to $\delta$, the substitution $4\mapsto 031$ violates the
 conditions of (ii) above. Other instances of the correspondence are
\[(2,1,1,2,2,2,1,1,2,2) \leftrightarrow (2,1,7,1,5)  \text{ \ and \ }
(1,2,2,1,2,1,1,2,2) \leftrightarrow (0,5,3,1,5)\,.\]

We  use these correspondences  to label basis elements of $\pppn$ by thin or almost-odd compositions instead of sparse subsets. Thus, given a thin composition $\tau$ of $n$ we let $Q_{\tau}:=Q_F$, where $F$ is the sparse subset of $[n-1]$ such that
$\tau_F=\tau$, and given an almost-odd composition $\gamma$ of $n$ we let $Q_{\gamma}:=Q_F$, where $F$ is the sparse subset of $[n-1]$ such that $\gamma_F=\gamma$, and similarly for the other bases.

\begin{example}\label{Ex:specialTO} Suppose $n$ is even. 
The almost-odd composition $(n)$ corresponds to 
the sparse subset   $\{1,3,5,\ldots,n-1\}$ and to the thin composition
$(\underbrace{2,2,\ldots,2}_{\frac{n}{2}})$. Thus,
 \begin{align*}
 \TO_{(n)}&=\TO_{(2,2,\ldots,2)}=\TO_{\{1,3,5,\ldots,n-1\}}=P_\emptyset\,,\\
  O_{(n)}&=O_{(2,2,\ldots,2)}=O_{\{1,3,5,\ldots,n-1\}}=\sum_{G\subseteq\{2,4,\ldots,n-2\}} P_G\,,\\
 Q_{(n)}&=Q_{(2,2,\ldots,2)}=Q_{\{1,3,5,\ldots,n-1\}}= P_{\{1,3,5,\ldots,n-1\}}\,.
\end{align*}
 If $n$ is odd,  the almost-odd composition $(0,n)$ corresponds to the
  sparse subset   $\{2,4,6,\ldots,n-1\}$ and to the thin composition
  $(1,\underbrace{2,\ldots,2}_{\frac{n-1}{2}})$. Thus,
 \begin{align*}
\TO_{(0,n)} &=\TO_{(1,2,\ldots,2)}=\TO_{\{2,4,6,\ldots,n-1\}}=P_\emptyset+P_{\{1\}}\,,\\   
O_{(0,n)}&=O_{(1,2,\ldots,2)}=O_{\{2,4,6,\ldots,n-1\}}=\sum_{G\subseteq\{1,3,\ldots,n-2\}} P_G\,,\\
 Q_{(0,n)}&=Q_{(1,2,\ldots,2)}=Q_{\{2,4,6,\ldots,n-1\}}= P_{\{2,4,6,\ldots,n-1\}}\,.\\
\end{align*}
\end{example}

Lemma~\ref{L:sparse-thin} allows us to rewrite~\eqref{E:OinQ} as follows: for any thin composition $\tau$ of $n$,
\begin{equation}\label{E:OinQtwo}
O_\tau=\!\!\!\!\sumsub{\rho \text{ thin}\\{\tau\leq\rho}}(-1)^{n-k(\rho)}Q_\rho\,.
\end{equation}

Lemma~\ref{L:sparse-almostodd} allows us to rewrite formulas~\eqref{E:XinQ} and~\eqref{E:TOinQ} as follows
(recall our convention~\eqref{E:parts} on the number of parts): for any pseudocomposition $\beta$ of $n$,
\begin{align}
\label{E:XinQtwo}
\varphi(X_\beta) &=2^{k(\beta)}\cdot\!\!\!\!\sumsub{\gamma \text{ almost-odd}\\ \beta\leq\gamma} (-1)^{\frac{n-k(\gamma)}{2}}Q_\gamma\,,\\
\intertext{and for any almost-odd composition $\gamma$ of $n$,}
\label{E:TOinQtwo}
\TO_\gamma &=\!\!\!\!\sumsub{\delta \text{ almost-odd}\\ \gamma\leq\delta} (-1)^{\frac{n-k(\delta)}{2}}Q_\delta\,.
\end{align}

\begin{corollary}\label{C:XontoTO} For any almost-odd composition $\gamma$ of $n$,
\begin{equation}\label{E:XontoTO}
\varphi(X_\gamma)=2^{k(\gamma)}\cdot \TO_\gamma\,.
\end{equation}
\end{corollary}\qed

\begin{remark}\label{R:hsiao-schocker}
Some of the definitions and results of Sections~\ref{S:bases} and~\ref{S:compositions} have counterparts in earlier work of Hsiao~\cite{Hsi} and Schocker~\cite{Sch}. These references do not deal with the peak algebra but with the {\em peak ideal}. Our study of the peak algebra is more general, although the underlying combinatorics is similar for both situations (almost-odd compositions versus odd compositions). See Remark~\ref{R:hsiao-schocker-two} for more
details.
\end{remark}

 \section{Chains of ideals of $\SolB$ and of $\pppn$}\label{S:chains}

For $n\geq 2$, consider the map $\pi_n:\pppn\to\ppp{n-2}$ defined by
\begin{equation}\label{E:defpi}
P_F\mapsto
\begin{cases}
-P_{F\setminus\{1\}-2} & \text{ if }1\in F\\
0 & \text{ if }2\in F\\
P_{F-2} & \text{ if neither $1$ nor $2$ belong to $F$}
\end{cases}
\end{equation}
for any sparse subset $F\subseteq [n-1]$. We let 
 $\pi_1$ and $\pi_0$ be the zero maps on $\ppp{1}$ and
 $\ppp{0}$, respectively. We often omit the
subindex $n$ from $\pi_n$. We know that $\pi:\pppn\to\ppp{n-2}$ is
a surjective morphism of algebras~\cite[Proposition 5.6]{ABN}.

We  describe $\pi$ on the other bases of the peak algebra (Definition~\ref{D:bases}).

\bp\label{P:pi-bases} Let $F$ be a sparse subset of $[n-1]$, $(a_1,\ldots,a_k)$ a thin composition of $n$, and $(b_0,b_1,\ldots,b_k)$ an almost-odd composition of $n$.
We have
\begin{align}
\label{E:pi-Q}
\pi(Q_F) & =\begin{cases} 
-Q_{F\setminus\{1\}-2} & \text{ if }1\in F\,,\\
                          0 & \text{ if }1\notin F\,;
             \end{cases}\\
\label{E:pi-O}
\pi(O_{(a_1,\ldots,a_k)}) & =\begin{cases} 
O_{(a_2,\ldots,a_k)} & \text{ if }a_1=2\,,\\
                          0 & \text{ if }a_1=1\,;
             \end{cases}\\
 \label{E:pi-TO}
 \pi(\TO_{(b_0,b_1,\ldots,b_k)}) & =\begin{cases} 
\TO_{(b_0-2,b_1,\ldots,b_k)} & \text{ if }b_0\geq 2\,,\\
                          0 & \text{ if }b_0=0\,.
             \end{cases}
\end{align}
\ep
\bpf By~\eqref{E:QinP}, $\pi(Q_F)=\sum_{F\subseteq G}\pi(P_G)$.
The only terms that contribute to this sum are those for which $2\notin G$.
These split in two classes: (i) those in which $1\in G$, and (ii) those in which $1,2\notin G$.
From~\eqref{E:defpi} we obtain
\[\pi(Q_F)=-\sum_{F\subseteq G,\,1\in G} P_{G\setminus\{1\}-2}+\sum_{F\subseteq G,\,1,2\notin G} P_{G-2}\,.\]
If $1\notin F$ there is a bijection from class (i) to class (ii) given by $G\mapsto G\setminus\{1\}$, and $\pi(Q_F)=0$. If $1\in F$ then class (ii) is empty and  class (i) is in bijection with the sparse subsets of $[n-3]$ which contain $F\setminus\{1\}-2$ via
$G\mapsto G\setminus\{1\}-2$; therefore, $\pi(Q_F)=-Q_{F\setminus\{1\}-2}$.

Let $\tau=(a_1,\ldots,a_k)$. Let $F$ be the sparse
subset of $[n-1]$ corresponding to $\tau$ as in Lemma~\ref{L:sparse-thin}, i.e.,
 $I(\tau)=[n-1]\setminus F$. If $a_1=1$ then $1\notin F$ and from~\eqref{E:OinQ}
 and~\eqref{E:pi-Q} we deduce $\pi(O_\tau)=0$. Assume $a_1=2$ and let
 $\hat{\tau}=(a_2,\ldots,a_k)$. Then $1\in F$ and $I(\hat{\tau})=I(\tau)-2=[n-3]\setminus(F\setminus\{1\}-2)$. We have
\begin{align*}
\pi(O_\tau) & \equal{\eqref{E:OinQ}} \sum_{G\subseteq F}(-1)^{\#G}\pi(Q_{G})\equal{\eqref{E:pi-Q}}
-\sum_{1\in G\subseteq F}(-1)^{\#G}Q_{G\setminus\{1\}-2}\\
&\ =\sum_{G'\subseteq F\setminus\{1\}-2}(-1)^{\#G'}Q_{G'}
\equal{\eqref{E:OinQ}}O_{\hat{\tau}}\,.
\end{align*}
The proof of~\eqref{E:pi-TO} is similar.
\epf

\begin{definition}\label{D:peakideals}  For each $j=0,\ldots,\ipartn$ let
$\pppn^j = \ker(\pi^{j+1}:\pppn\to\ppp{n-2j-2})$. 
\end{definition}

Since $\pi$ is a morphism of algebras, these subspaces form a chain of ideals
\[\pppn^0 \subseteq \pppn^1  \subseteq
\cdots \subseteq \pppn^{\ipartn}= \pppn\,.\]
In particular, $\ker(\pi)=\pppn^0$ is the {\em peak ideal}~\cite[Theorem 5.7]{ABN}.
This space has  a linear basis consisting of sums of permutations with a common set of {\em interior} peaks~\cite[Definition 5.5]{ABN}.

{}From Proposition~\ref{P:pi-bases} we deduce the following explicit
description of the ideals $\pppnj$.

\begin{corollary}\label{C:peakideals} Let $j=0,\ldots,\ipartn$ let $\pppnj$.
The ideal $\pppnj$ is linearly spanned by any of the sets consisting of:
\begin{itemize}
\item[(a)] The elements $Q_F$ as $F$ runs over the sparse subsets of $[n-1]$
which do not contain $\{1,3,\ldots,2j+1\}$.
\item[(b)] The elements  $O_{\alpha}$ as $\alpha=(a_1,\ldots,a_k)$ runs over 
those thin compositions of $n$ such that either $k\leq j$ or else
there is at least one index  $i\leq j+1$ with $a_i =1$.
\item[(c)] The elements  $\TO_{\beta}$ as $\beta=(b_0,b_1,\ldots,b_k)$ runs over 
those almost-odd compositions of $n$ such that $b_0\leq 2j$.
\end{itemize}
\end{corollary}\qed

The almost-odd compositions of $n$ that do not satisfy  condition (c)  are in bijection with the almost-odd compositions of
 $n-2j-2$ via $(b_0,b_1,\ldots,b_k)\mapsto (b_0-2j-2,b_1,\ldots,b_k)$.
 Therefore,
\[\dim\pppnj=\begin{cases} f_n-f_{n-2j-2} & \text{ if }j<\ipartn\,,\\
f_n & \text{ if }j=\ipartn\,.\end{cases}\]
The  sparse subsets of $[n-1]$ that do not satisfy condition (a) are those of the form $\{1,3,\ldots,2j+1\}\cup G$, where $G$ is a sparse subset of $\{2j+3,\ldots,n-1\}$.
The thin compositions $\alpha=(a_1,\ldots,a_k)$ that do not satisfy 
condition (b) are those for which $k\geq j+1$ and
$a_1=\ldots=a_{j+1}=2$. 

\begin{remark}\label{R:hsiao-schocker-two} 
Specializing $j=0$ in the preceding remarks   we obtain that the peak ideal $\pppn^0$  is linearly spanned by
\begin{itemize}
\item[(a)] The elements $Q_F$ as $F$ runs over the sparse subsets of $[n-1]$
which do not contain $1$.
\item[(b)] The elements  $O_{\alpha}$ as $\alpha=(a_1,\ldots,a_k)$ runs over 
those thin compositions of $n$ such that  $a_1 =1$.
\item[(c)] The elements  $\TO_{\beta}$ as $\beta=(b_0,b_1,\ldots,b_k)$ runs over 
those almost-odd compositions of $n$ such that $b_0=0$.
\end{itemize}
The dimension of the peak ideal is $\dim\pppn^0=f_n-f_{n-2}=f_{n-1}$.

The peak ideal  is the object studied in~\cite{BHT,Nym,Sch} (and in dual form in~\cite{Hsi}). The bases $Q$ and $\TO$ specialized as in (a) and (c) above are the bases $\Gamma$ and $\tilde{\Xi}$ of~\cite[Section 3]{Sch}. The basis $O$ appears to be new, even after specialization.
\end{remark}

\medskip

For $n\geq 1$, consider the map $\beta_n:\SolB\to\Sol{B_{n-1}}$  defined by
\begin{equation}\label{E:defbeta}
X_{(b_0,b_1,\ldots,b_k)} \mapsto 
\begin{cases} X_{(b_0-1,b_1,\ldots,b_k)} & \text{ if }b_0\neq 0\,,\\
                          0 & \text{ if }b_0=0\,.
             \end{cases}
\end{equation} 
We let $\beta_0$ be the zero map on $\Sol{B_0}$. We often omit the subindex $n$ from $\beta_n$.

\bd For each $i=0,\ldots,n$, let $\I_n^i=\ker(\beta^{i+1}:\SolB\to\Sol{B_{n-1}})$. 
\ed

We know that $\beta$ is a surjective morphism of algebras~\cite[Proposition 5.2]{ABN}. Therefore, these subspaces form a chain of ideals
\[\I_n^0 \subseteq \I_n^1 \subseteq \I_n^2 \subseteq
\cdots \subseteq \I_n^n = \SolB\,.\]

{}From~\eqref{E:defbeta} we deduce that the ideal $\I_n^i$ is  linearly spanned by
the elements $X_\beta$ as $\beta=(b_0,b_1,\ldots,b_k)$ runs over 
those pseudocompositions of $n$ such that $b_0 \leq i$.  

Under the canonical map $\varphi:\SolB\to\SolA$, the ideals $\I_n^0$ and $\I_n^1$
both map onto the peak ideal $\pppn^0$~\cite[Theorem 5.9]{ABN}.
Furthermore, $\I_n^n =  \SolB$ maps onto the
peak algebra $\pppn$~\cite[Theorem 4.2]{ABN}. These results generalize as follows.

\begin{proposition}\label{P:ideals}
For each $i=0,\ldots,n$,
\[\varphi(\I_n^{i})=\pppn^{\ipart{i}}\,.\]
\end{proposition}
\begin{proof} Let $j=0,\ldots,\ipartn$.
Let $\gamma=(c_0,c_1,\ldots,c_k)$ be an almost-odd composition of $n$ such that $c_0\leq 2j$. Then $\varphi(X_\gamma)=2^{k(\gamma)}\cdot\TO_{\gamma}$ by~\eqref{E:XontoTO}. In addition, $X_\gamma\in \I_n^{2j}$, and by Corollary~\ref{C:peakideals}, these elements $\TO_\gamma$ span $\pppnj$.
Therefore, $\pppnj\subseteq \varphi(\I_n^{2j})$. 

On the other hand,  the commutativity of the diagram~\cite[Proposition 5.6]{ABN}
\[\xymatrix{ {\SolB}\ar[r]^{\beta^2}\ar[d]_{\varphi} &
{\Sol{B_{n-2}}}\ar[d]^{\varphi}\\
{\pppn}\ar[r]_{\pi}&{\ppp{n-2}} }\]
implies that $\I_n^{2j+1}=\ker(\beta^{2j+2})$ maps under $\varphi$ to $\ker(\pi^{j+1})=\pppnj$.

Thus $\pppnj\subseteq \varphi(\I_n^{2j})\subseteq \varphi(\I_n^{2j+1})\subseteq\pppnj$ and the result follows.
\end{proof}

The situation may be illustrated as follows:

 \begin{center}\begin{picture}(200,60)
\put(0,40){$\underbrace{\I_n^0 \subseteq \I_n^1} \subseteq 
 \underbrace{\I_n^2 \subseteq \I_n^3}\subseteq
\cdots \subseteq \I_n^n=\SolB$}
\put(17.5,20){$\downarrow$}\put(17.5,18){$\downarrow$}
\put(74,20){$\downarrow$}\put(74,18){$\downarrow$}
\put(184.8,20){$\downarrow$}\put(184.8,18){$\downarrow$}\put(186.3,25){$\mid$}
 \put(15,0){$\pppn^0 \ \ \ \subseteq\ \ \ 
 \pppn^1\ \ \ \subseteq
\cdots  \subseteq  \pppn^{\ipartn}=\pppn$}
\end{picture}\end{center}

\section{The radical of the peak algebra}\label{S:radical}

Let $A$ be an Artinian ring (e.g., a finite dimensional algebra over a field).
The {\em (Jacobson) radical} $\rad(A)$ may be defined in any of the following ways~\cite[Theorem 4.12, Exercise 11 in Section 4]{Lam}:
\begin{itemize}
\item[(R1)] $\rad(A)$ is the largest nilpotent ideal of $A$;
\item[(R2)] $\rad(A)$ is the smallest ideal of $A$ such that the corresponding quotient is semisimple. 
\end{itemize}
Thus, $\rad(A)$ is a nilpotent ideal and an ideal $N$ is nilpotent if and only if $N\subseteq \rad(A)$; $A/\rad(A)$ is semisimple  and an ideal $I$ is such that $A/I$ is semisimple if and only if $I\supseteq \rad(A)$.

\begin{lemma}\label{L:radical}
Let $A$ be an Artinian ring and $f:A\to B$ a surjective morphism of rings. Then
\[f(\rad(A))=\rad(B)\,.\]
\end{lemma}
\begin{proof} Since $f$ is surjective, $f(\rad(A))$ is an ideal of $B$. Since $\rad(A)$
is nilpotent, so is $f(\rad(A))$. Hence, by (R1), $f(\rad(A))\subseteq \rad(B)$. 
On the other hand, $f$ induces an isomorphism of rings
\[ A/f^{-1}(f(\rad(A)))\cong B/f(\rad(A))\,.\]
 Since $f^{-1}(f(\rad(A)))\supseteq \rad(A)$, the quotient is semisimple, by (R2) applied to $A$. Hence,  by (R2) applied to $B$, $f(\rad(A))\supseteq \rad(B)$. 
\end{proof}

We apply the lemma to derive an explicit description of the radical of the peak algebra from the known description of the radical of the descent algebra of type B. Solomon described the radical of the descent algebra of an arbitrary finite
Coxeter group~\cite[Theorem 3]{Sol}.  For the  descent algebra of type B, his result specializes as follows (see also~\cite[Corollary 2.13]{B92}).
   Given a pseudocomposition $\beta=(b_0,b_1,\ldots,b_k)$ of $n$ and
a permutation $s$ of $[k]$, let
\begin{equation}\label{E:beta-s}
\beta^s:=(b_0,b_{s(1)},\ldots,b_{s(k)})\,.
\end{equation}
The radical $\rad(\SolB)$ is linearly spanned by the
elements
\begin{equation}\label{E:radX}
X_{\beta}-X_{\beta^s}
\end{equation}
as $\beta$ runs over all pseudocompositions of $n$ and $s$ over all permutations of $[k(\beta)]$. It follows that the dimension of the
maximal semisimple quotient of $\SolB$ is 
\begin{equation}\label{E:codimradB}
\codim \rad(\SolB)=p(0)+p(1)+\cdots+p(n)\,, 
\end{equation}
where $p(n)$ is the number of partitions of $n$.

\begin{theorem}\label{T:radical}
The radical $\rad(\pppn)$ is linearly spanned by the elements in either (a), (b), or (c):
\[ \text{(a)}\ \ \TO_{\gamma}-\TO_{\gamma^t}\,,\ \ 
\text{(b)}\ \ Q_{\gamma}-Q_{\gamma^t}\,,\ \ 
\text{(c)}\ \ O_{\gamma}-O_{\gamma^t}\,.\ \ \]
In each case, $\gamma$ runs over all almost-odd compositions of $n$ and $t$ over all permutations of $[k(\gamma)]$.
\end{theorem}
\begin{proof} Let $J_a$, $J_b$, and $J_c$ be the span of the elements in $(a)$,  
$(b)$, and $(c)$ respectively.

Consider the canonical morphism $\varphi:\SolB\to\SolA$. Its image is $\pppn$~\cite[Theorem 4.2]{ABN}.
According to Lemma~\ref{L:radical} and~\eqref{E:radX}, $\rad(\pppn)$ is spanned by the elements
\[\varphi(X_{\beta})-\varphi(X_{\beta^s}) \]
with $\beta$ and $s$ as in~\eqref{E:beta-s}. 

Given an almost-odd composition $\gamma$ and a permutation $t$ of $[k(\gamma)]$,
~\eqref{E:XontoTO} gives
\[\varphi(X_\gamma)-\varphi(X_{\gamma^t})=2^{k(\gamma)}\cdot (\TO_\gamma-\TO_{\gamma^t})\,.\]
This shows that $J_a\subseteq \rad(\pppn)$.

Fix $\beta$ and $s$ as in~\eqref{E:beta-s}. Given a pseudocomposition $\gamma\geq\beta$, write $\gamma=\gamma_0\gamma_1\cdots\gamma_k$ (concatenation of compositions), with $\gamma_0$ a pseudocomposition of $b_0$ and $\gamma_i$ an ordinary composition of $b_i$ for $i=1,\ldots,k$. Define
\[\gamma^s:=\gamma_0\gamma_{s(1)}\cdots\gamma_{s(k)}\,.\]
This extends definition~\eqref{E:beta-s}. Note that $\gamma^s\geq\beta^s$, and  if $\gamma$ is almost-odd then so is $\gamma^s$. Therefore,
the map $\gamma\mapsto\gamma^s$ is a bijection from the almost-odd compositions
$\gamma\geq\beta$ to the almost-odd compositions $\gamma'\geq\beta^s$ (the inverse is $\gamma'\mapsto(\gamma')^{s^{-1}}$). Note also that $k(\gamma)=k(\gamma^s)$. Together with~\eqref{E:XinQtwo}
this gives
\[\varphi(X_{\beta})-\varphi(X_{\beta^s})=2^{k(\beta)}\cdot\!\!\!\!\sumsub{\gamma \text{ almost-odd}\\ \beta\leq\gamma} (-1)^{\frac{n-k(\gamma)}{2}}(Q_\gamma-Q_{\gamma^s})\,.\]
Note that each $\gamma^s=\gamma^t$ for a certain permutation $t$ of $[k(\gamma)]$.
This shows that $\rad(\pppn)\subseteq J_b$.

The bijection of the preceding paragraph may also be used in conjunction with~\eqref{E:TOinQtwo} to give
\[\TO_\gamma -\TO_{\gamma^s}=\!\!\!\!\sumsub{\delta \text{ almost-odd}\\ \gamma\leq\delta} (-1)^{\frac{n-k(\delta)}{2}}(Q_\delta-Q_{\delta^s})\,.\]
M\"obius inversion then shows that $J_b\subseteq J_a$.
Thus $J_a=\rad(\pppn)=J_b$.

Lastly, we deal with $J_c$. Recall the bijection $\gamma\mapsto\tau_\gamma$ between almost-odd compositions and thin compositions of Lemma~\ref{L:thin-almostodd}. Consider~\eqref{E:OinQtwo}. When written in terms of
almost-odd compositions, this equation says that
\[O_\gamma=\sum_{\delta}(-1)^{\frac{n-k(\delta)}{2}}Q_\delta\]
the sum being over those almost-odd compositions $\delta$ such that $\tau_\gamma\leq\tau_\delta$. Let $t$ be a permutation of $[k(\gamma)]$.
The map $\delta\mapsto\delta^t$ restricts to a bijection between the almost-odd
compositions $\delta$ such that $\tau_\gamma\leq\tau_\delta$ and the almost-odd
compositions $\delta'$ such that $\tau_{\gamma^t}\leq\tau_{\delta'}$. This is so because the restriction on the admissible refinements described in item (ii) of Lemma~\ref{L:thin-almostodd} only depends on the individual parts of $\gamma$, and not on their relative position.
Therefore,
\[O_\gamma-O_{\gamma^t}=\sum_{\delta}(-1)^{\frac{n-k(\delta)}{2}}(Q_\delta-Q_{\delta^t})\,.\]
Together with M\"obius inversion this shows that $J_c=J_b$.
\end{proof}

A partition of $n$ is an ordinary composition $\lambda=(\ell_1,\ell_2,\ldots,\ell_k)$ of $n$ such that $\ell_1\geq\ell_2\geq\ldots\geq\ell_k$.
 We say that  $\lambda$ is odd if each $\ell_i$ is odd, and almost-odd if at most one $\ell_i$ is even.

\begin{corollary}\label{C:radical} The dimension of the maximal semisimple quotient
of $\pppn$ is the number of almost-odd partitions of $n$.
\end{corollary}\qed

An almost-odd partition of $n$ may be viewed as an odd partition of $m$ for some $m\leq n$ such that $n-m$ is even. Therefore, the dimension of the maximal semisimple quotient of $\pppn$ is
\begin{equation}\label{E:codimradP}
\codim \rad(\pppn)=p_o(n)+p_o(n-2)+p_o(n-4)+\cdots+p_o(n-2\ipartn)\,,
\end{equation}
 where $p_o(n)$ is the number of odd partitions of $n$.
The number of almost-odd partitions of $n$ is, for $n\geq 0$,
\[1,1,2,3,4,6,8,11,14,19,\ldots.\]
For more information on this sequence, see~\cite[A038348]{SlSeek}.

The partial sums of~\eqref{E:codimradB} and~\eqref{E:codimradP} are
the codimensions of the radicals of the ideals of Section~\ref{S:chains}.
\begin{corollary}\label{C:radical-ideals} For any $i=0,\ldots,n$,
\[\dim \frac{\I_n^i}{\rad(\SolB)\cap\I_n^i}=p(n)+p(n-1)+\cdots+p(n-i)\]
 and for any $j=0,\ldots,\ipartn$,
\[\dim \frac{\pppn^j}{\rad(\pppn)\cap\pppn^j}=p_o(n)+p_o(n-2)+\cdots+p_o(n-2j)\,.\]
\end{corollary}
\begin{proof} The first equality follows directly from~\eqref{E:radX} and the
definition of the ideals $\I_n^i$. The second follows from Theorem~\ref{T:radical}
and item (c) in Corollary~\ref{C:peakideals}.
\end{proof}

In particular,  the codimension  of the radical of the peak ideal $\pppn^0$ is the number of odd partitions of $n$. This result is due to  Schocker~\cite[Corollary 10.3]{Sch}. (In this reference, $\pppn^0$ is viewed as a non-unital algebra, but 
this leads to the same answer, since the radical of an ideal of a ring coincides with the intersection of the ideal with the radical of the ring~\cite[Exercise 7 in Section 4]{Lam}.)

\medskip

The radical may also be described in terms of thin compositions in either of the three bases, by transporting the action~\eqref{E:beta-s} of permutations on almost-odd compositions to an action on thin compositions via the bijection of Lemma~\ref{L:thin-almostodd}. We describe the result.

Given a thin composition $\tau$, consider the unique way of writing it as the concatenation of compositions $\tau=\tau_0\tau_1\cdots\tau_h$ in which $\tau_0$ is of the form $(2,2,\ldots,2)$ ($\tau_0$ may be empty), and for each $i>0$  $\tau_i$ is of the form $(1,2,\ldots,2)$. 
For instance, if $\tau=(2,1,1,2,2,2,1,1,2,2)$ then $\tau_0=(2)$, $\tau_1=(1)$, $\tau_2=(1,2,2,2)$, $\tau_3=(1)$, $\tau_4=(1,2,2)$. Let $h:=h(\tau)$. (If $\tau=\tau_\gamma$ then $h(\tau)=k(\gamma)$.) Given a permutation $t$ of $[h(\tau)]$ we let $\tau^t:=\tau_0\tau_{t(1)}\cdots\tau_{t(h)}$.

\begin{proposition}\label{P:radical}
The radical $\rad(\pppn)$ is linearly spanned by the elements in either (a), (b), or (c):
\[ \text{(a)}\ \ \TO_{\tau}-\TO_{\tau^t}\,,\ \ 
\text{(b)}\ \ Q_{\tau}-Q_{\tau^t}\,,\ \ 
\text{(c)}\ \ O_{\tau}-O_{\tau^t}\,.\ \ \]
In each case, $\tau$ runs over all thin compositions of $n$ and $t$ over all permutations of $[h(\tau)]$.
\end{proposition}
\qed

\begin{remark}\label{R:radical}
We point out that the radicals of the peak algebra and the descent algebra of type A are related by
\[\rad(\pppn)=\rad(\SolA)\cap\pppn\,.\]
First, for any extension of algebras $A\subseteq B$, we have that $\rad(B)\cap A$ is a nilpotent ideal of $A$, so $\rad(B)\cap A\subseteq\rad(A)$ by (R1). The reverse inclusion does not always hold, but it does if $B/\rad(B)$ is commutative.  Indeed, a commutative semisimple algebra does not contain nilpotent elements, and since  $A/(\rad(B)\cap A)\inc B/\rad(B)$, $A/(\rad(B)\cap A)$ does not contain nilpotent elements. Hence $A/(\rad(B)\cap A)$ is semisimple by (R1), and then
$\rad(A)\subseteq \rad(B)\cap A$ by (R2).
These considerations apply in our situation ($A=\pppn$, $B=\SolA$), since it is known that $\SolA/\rad(\SolA)$ is commutative~\cite[Theorem 3]{Sol} (this quotient is isomorphic to the representation ring of $S_n$).
\end{remark}

\section{The convolution product}\label{S:convolution}

The convolution product of permutations is due to Malvenuto and Reutenauer~\cite{MalReu}. It may also be defined for signed permutations.
We review the relevant notions below, for more details see~\cite[Section 8]{ABN}.

Consider the spaces
\[\field B:=\bigoplus_{n\geq 0}\field B_n \text{ \ and \ }
\field S:=\bigoplus_{n\geq 0}\field S_n\,.\]
On the space $\field S$ there is defined the {\em external} or {\em convolution} product
\begin{equation}\label{E:convolution}
\sigma\ast \tau:=\sum_{\xi\in\Sh(p,q)}\xi\cdot(\sigma\times \tau)\,.
\end{equation}
Here $\sigma\in S_p$ and $\tau\in S_q$ are  permutations,
\[\Sh(p,q)=\{\xi\in S_{p+q} \mid \xi(1)<\cdots<\xi(p),\ \xi(p+1)<\cdots<\xi(p+q)\}\]
 is the set of $(p,q)$-{\em shuffles}, and $\sigma\times \tau\in S_{p+q}$ is defined by 
 \[(\sigma\times \tau)(i)=\begin{cases} \sigma(i)
 & \text{ if }1\leq i\leq p\,,\\
 p+\tau(i-p) & \text{ if }p+1\leq i\leq p+q\,. \end{cases}\]
 The convolution product turns the space $\field S$ into a graded algebra.
 
Similar formulas define the convolution product on $\field B$, and the canonical map $\varphi:\field B\to \field S$ preserves this structure. 

Consider now the spaces
\[\Sol{B}:=\bigoplus_{n\geq 0}\SolB\,,\quad
\ppp{}:=\bigoplus_{n\geq 0}\pppn\,,\quad
\I^0:=\bigoplus_{n\geq 0}\I_n^0\,,\quad
\ppp{}^0:=\bigoplus_{n\geq 0}\pppn^0\,.\]

Under the convolution product of $\field B$, $\I^0$ is a graded subalgebra of $\field B$ and $\Sol{B}$ is a graded right $\I^0$-submodule of $\field B$. Similarly, $\ppp{}^0$ is a graded subalgebra of $\field S$ and $\ppp{}$ is a graded right $\ppp{}^0$-submodule of $\field S$, and the map $\varphi$ preserves each of these structures. The situation may be schematized by
\[\xymatrix@C=0pc{ {\I^{0}}\ar@{->>}[d]_{\varphi} &\subseteq &
{\Sol{B}}\ar@{->>}[d]_{\varphi} &\subseteq& {\field B}\ar@{->>}[d]^{\varphi}\\
 {\ppp{}^0} &\subseteq& {\ppp{}}& \subseteq & {\field S} }\]

For any pseudocomposition $(b_0,b_1,\ldots,b_k)$  we have 
\begin{equation}\label{E:Xfree}
X_{(b_0,b_1,\ldots,b_k)}=X_{(b_0)}\ast X_{(0,b_1)}\ast\cdots\ast X_{(0,b_k)}\,.
\end{equation}
In other words, the basis $X$ of $\Sol{B}$ is {\em multiplicative} with respect to convolution. 

We deduce that all three bases $Q$, $O$, and $\TO$ of $\ppp{}$ are
multiplicative.

\begin{proposition}\label{P:peakfree} For any almost-odd composition $(b_0,b_1,\ldots,b_k)$,
\begin{align}\label{E:TOfree}
\TO_{(b_0,b_1,\ldots,b_k)} &=\TO_{(b_0)}\ast \TO_{(0,b_1)}\ast\cdots\ast \TO_{(0,b_k)}\,,\\
\label{E:Qfree}
Q_{(b_0,b_1,\ldots,b_k)} &=Q_{(b_0)}\ast Q_{(0,b_1)}\ast\cdots\ast Q_{(0,b_k)}\,,\\
\label{E:Ofree}
O_{(b_0,b_1,\ldots,b_k)} &=O_{(b_0)}\ast O_{(0,b_1)}\ast\cdots\ast O_{(0,b_k)}\,.
\end{align}
\end{proposition}
\begin{proof} Formula~\eqref{E:TOfree} follows at once from~\eqref{E:Xfree} by applying the canonical map $\varphi$, in view of~\eqref{E:XontoTO} and the fact that $\varphi$ preserves the convolution product.

Let $\beta:=(b_0,b_1,\ldots,b_k)$. From~\eqref{E:TOinQtwo} we obtain
\[\TO_{(b_0)}\ast \TO_{(0,b_1)}\ast\cdots\ast \TO_{(0,b_k)}=\sumsub{(b_0)\leq\delta_0\\(0,b_i)\leq\delta_i}(-1)^{\frac{n-\sum_{i=0}^k k(\delta_i)}{2}}Q_{\delta_0}\ast Q_{\delta_1}\ast\cdots\ast Q_{\delta_k}\,.\]
The sum is over almost-odd compositions $\delta_0$ and $\delta_i$ as indicated. For $i>0$, any such $\delta_i$ is of the form $(0,\alpha_i)$ with $\alpha_i$ an (ordinary) odd composition of $b_i$. Note $k(\alpha_i)=k(\delta_i)$~\eqref{E:parts}.
The concatenation $\delta:=\delta_0\alpha_1\ldots\alpha_k$ is then an almost-odd composition of $n$ with $k(\delta)=\sum_{i=0}^k k(\delta_i)$ and $\delta\geq\beta$. Any almost-odd composition $\delta\geq\beta$ is of this form
for a unique sequence $\delta_i$. Therefore, the right hand side may be written as
\[\sum_{\beta\leq\delta}(-1)^{\frac{n- k(\delta)}{2}}Q_{\delta_0}\ast Q_{\delta_1}\ast\cdots\ast Q_{\delta_k}\,.\]
On the other hand, by\eqref{E:TOinQtwo} and~\eqref{E:TOfree}, the left hand side equals
\[\sum_{\beta\leq\delta}(-1)^{\frac{n- k(\delta)}{2}}Q_{\delta}\,.\]
By M\"obius inversion we deduce $Q_{\delta_0}\ast Q_{\delta_1}\ast\cdots\ast Q_{\delta_k}=Q_\delta$ for each $\delta$, which gives~\eqref{E:Qfree}.

Formula~\eqref{E:Ofree} may be deduced similarly from~\eqref{E:OinQtwo} and~\eqref{E:Qfree}. The argument now involves the partial order on almost-odd compositions  corresponding to refinement of thin compositions. Let $\beta$, $\delta$, and $\delta_i$ be as above. The key observation is that $\beta\leq\delta$ in this partial order if and only if $(b_0)\leq\delta_0$ and $(0,b_i)\leq\delta_i$ ($i>0$) in the same partial order. This is guaranteed by
item (ii) of Lemma~\ref{L:thin-almostodd}.
\end{proof}

Equation~\eqref{E:Xfree} implies
\begin{align*}\label{E:prodSolB}
X_{(b_0,b_1,\ldots,b_k)}\ast X_{(0,c_1,\ldots,c_h)}=X_{(b_0,b_1,\ldots,b_k,c_1,\ldots,c_h)}\,,
\intertext{and in particular,}
X_{(0,a_1,\ldots,a_k)}\ast X_{(0,c_1,\ldots,c_h)}=X_{(0,a_1,\ldots,a_k,c_1,\ldots,c_h)}\,.
\end{align*}
The second equation says that $\I^0$ is a free algebra, with one generator of degree $n$ for each $n$ (the element $X_{(0,n)}$). The first equation says that
$\Sol{B}$ is a free $\I^0$-module, with one generator of degree $n$ for each
$n$ (the element $X_{(n)}$). The latter fact is reflected in the following relation between the Hilbert series of these graded vector spaces:
\[\frac{\Sol{B}(t)}{\I^0(t)}=\frac{1}{1-t}\,.\]

Similarly, Proposition~\ref{P:peakfree} implies that $\ppp{}^0$ is a free algebra with one generator of degree $n$ for each odd $n$ (a result known from~\cite{BMSW,Hsi,Sch}) and also that
 $\ppp{}$ is a free $\ppp{}^0$-module, with one generator of degree $n$ for each even $n$.
Correspondingly, the Hilbert series of these graded vector spaces are related by
\[\frac{\ppp{}(t)}{\ppp{}^0(t)}=\frac{1}{1-t^2}\,.\]

\section{Eulerian idempotents and a  basis of semiidempotents}\label{S:basis}

\subsection{Peak analogs of the first Eulerian idempotent}\label{S:first}
As in~\cite[Section 6]{ABN},
we consider certain elements of the group algebras of $B_n$ and $S_n$ obtained by grouping permutations according to their number of descents, peaks, interior descents, or interior peaks, respectively. More precisely, we let
\begin{align*}
y_j \  & :=\sum\{\sigma \in B_n \mid \#\Des(\sigma)=j\} \text{ \ \ for \ }j=0,\ldots,n\,;\\
y^0_j\  & :=\sum\{\sigma \in B_n\mid \#(\Des(\sigma)\setminus \{0\})=j-1\}
 \text{ \ \ for \ }j=1,\ldots,n\,;\\
 p_j\ &:= \sum\{\sigma \in S_n \mid \#\Peak(\sigma)=j\}  \text{ \ \ for \ } j=0,\ldots,\ipartn\,;\\
p^0_j\ & :=\sum\{\sigma \in S_n \mid \#(\Peak(\sigma)\setminus \{1\})=j-1\}
 \text{ \ \ for \ } j=1,\ldots,\ipart{n+1}\,.
\end{align*}
We have that $y_j\in\SolB$, $y_j^0\in\I_n^0$, $p_j\in\pppn$, and $p_j^0\in\pppn^0$.  The canonical map $\varphi:\SolB\to\pppn$ satisfies~\cite[Propositions 6.2, 6.4]{ABN}
\begin{align} \label{E:yinp}
\varphi(y_j) &= \sum_{i=0}^{\min(j,n-j)} 2^{2i} \binom{n-2i}{j-i}\cdot p_i
\text{ \ \ for \ } j=0,\ldots,n\,;\\
\label{E:y0inp0}
\varphi(y_j^0)  &= \sum_{i=1}^{\min(j,n+1-j)} 2^{2i-1}\binom{n-2i+1}{j-i}\cdot p^0_i \text{ \ \ for \ }j=1,\ldots,n\,.
\end{align}

For each $n\in\Z^+$, let $n!!:=n(n-2)(n-4)\cdots$ (the last term
in the product is $2$ if $n$ is even and $1$ if $n$ is odd).
Set also $0!!=(-1)!!=1$. Note that 
\begin{equation}\label{E:dfactorials}
(2n)!!=2^nn! \text{ \ and \ } (2n+1)!!=\frac{(2n+1)!}{2^nn!}\,.
\end{equation}

Consider the following elements:
\begin{align}
\label{E:def-e} e_{(n)} &=\sum_{j=0}^n (-1)^j\frac{(2j-1)!!(2n-2j-1)!!}{(2n)!!}y_j\in\SolB\\
\label{E:def-e0} e_{(0,n)} &=\sum_{j=1}^n (-1)^{j-1}\frac{(j-1)!(n-j)!}{n!}y^0_j\in \I_n^0\\
\label{E:def-rho} \rho_{(n)} &=\sum_{j=0}^{\ipartn}(-1)^j\frac{(2j-1)!!(n-2j-1)!!}{n!!}p_j\in\pppn\\
\label{E:def-rho0} \rho_{(0,n)} &=\sum_{j=1}^{\ipart{n+1}}(-1)^{j-1}\frac{(2j-2)!!(n-2j)!!}{n!!}\pint_j\in\pppn^0
\end{align}

These elements are analogous to a certain element of $\SolA$ known as the {\em first Eulerian idempotent}.
The elements $e_{(n)}$ and $e_{(0,n)}$ appear in work of Bergeron and Bergeron~\cite{BB92a,BB92b,B92}, where they are denoted $I_\emptyset$ and $I_{(n)}$, respectively.  According to~\cite[Theorems 2.1, 2.2]{B92}, $e_{(n)}$ and $\frac{1}{2}e_{(0,n)}$ are orthogonal idempotents. For odd $n$, the element $\rho_{(0,n)}$ is known to be idempotent from work of Schocker~\cite[Section 7]{Sch}. Below we deduce this fact, as well as the idempotency of $\rho_{(n)}$ for even $n$, from that of $e_{(n)}$ and $e_{(0,n)}$.

\begin{lemma}\label{L:binomials}
For any $i=1,\ldots,\ipart{n+1}$, 
\begin{equation}\label{E:binomial1}
\sum_{j=i}^{n-i+1} (-1)^{j-i}
\frac{(j-1)!(n-j)!}{(j-i)!(n-i-j+1)!} =\begin{cases}
\frac{(2i-2)!!(n-1)!!}{2^{2i-2} (n-2i+1)!!} & \text{ if $n$ is odd,}\\
0 & \text{ if $n$ is even.}
\end{cases}
\end{equation}
For any $i=0,\ldots,\ipartn$, 
\begin{equation}\label{E:binomial2}
\sum_{j=i}^{n-i} (-1)^{j-i} \frac{(2j-1)!(2n-2j-1)!}{(j-i)!(j-1)!(n-i-j)!(n-j-1)!} =\begin{cases}
\frac{(2i-1)!!(n-1)!!}{2^{2i-2n+2}(n-2i)!!} & \text{ if $n$ is even,}\\
0 & \text{ if $n$ is odd.}
\end{cases}
\end{equation}
\end{lemma}
\begin{proof} Start from the equality
\[\frac{1}{(1+x)^i}\cdot\frac{1}{(1-x)^i}=\frac{1}{(1-x^2)^i}\,.\]
Expanding with the binomial theorem  gives
\[\sum_{r=0}^\infty (-1)^r\binom{i+r-1}{r}x^r\cdot \sum_{s=0}^\infty \binom{i+s-1}{s}x^s = \sum_{t=0}^\infty \binom{i+t-1}{t}x^{2t}\,.\]
Equating coefficients of $x^m$ we obtain
\[\sum_{r=0}^m (-1)^r\binom{i+r-1}{r} \binom{i+m-r-1}{m-r} = \begin{cases} \binom{i+m/2-1}{m/2} & \text{ if $m$ is even,}\\
0 & \text{ if $m$ is odd.}\end{cases}\]
Letting $n=m+2i-1$ and $j=r+i$ this equality becomes
\[\sum_{j=i}^{n-i+1} (-1)^{j-i}
\binom{j-1}{j-i}\binom{n-j}{n-i-j+1} =\begin{cases}
\binom{n/2-1/2}{n/2-i+1/2} & \text{ if $n$ is odd,}\\
0 & \text{ if $n$ is even.}
\end{cases}\]
Equation~\eqref{E:binomial1} follows by noting that $\binom{n/2-1/2}{n/2-i+1/2}(i-1)!^2=\frac{(2i-2)!!(n-1)!!}{2^{2i-2} (n-2i+1)!!}$ for odd $n$.

Equation~\eqref{E:binomial2} can be deduced similarly, starting from
\[\frac{1}{(1+x)^{i+\frac{1}{2}}}\cdot\frac{1}{(1-x)^{i+\frac{1}{2}}}=\frac{1}{(1-x^2)^{i+\frac{1}{2}}}\,.\]
\end{proof}

\begin{theorem}\label{T:etorho} Let $n$ be a positive integer. Then
\begin{align}
\label{E:etorho}
\varphi(e_{(n)}) &=\begin{cases}
\rho_{(n)} & \text{ if $n$ is even,}\\ 
0 & \text{ if $n$ is odd;}
\end{cases}\\
\label{E:e0torho0}
\varphi(e_{(0,n)}) &=\begin{cases}
0 & \text{ if $n$ is even,}\\ 
2\rho_{(0,n)} & \text{ if $n$ is odd.}
\end{cases}
\end{align}
In particular, $\rho_{(n)}$ is idempotent for each even $n$, and $\rho_{(0,n)}$
is idempotent for each odd $n$.
\end{theorem}
\begin{proof} According to~\eqref{E:yinp} and~\eqref{E:def-e},
\begin{align*}
\varphi(e_{(n)}) &= \sum_{j=0}^n (-1)^j
\frac{(2j-1)!!(2n-2j-1)!!}{(2n)!!} \sum_{i=0}^{\rm{min}(j,n-j)} 2^{2i}
{{n-2i}\choose{j-i}}\cdot p_i\\
&= \sum_{i=0}^{\lfloor \frac{n}{2} \rfloor} 2^{2i} \sum_{j=i}^{n-i}
(-1)^j \frac{(2j-1)!!(2n-2j-1)!!}{(2n)!!}{{n-2i}\choose{j-i}}\cdot p_i\\
& \equal{\eqref{E:dfactorials}} \sum_{i=0}^{\lfloor \frac{n}{2} \rfloor}
\frac{(n-2i)!}{n!} 2^{2i-2n+2}\sum_{j=i}^{n-i} (-1)^j \frac{(2j-1)!(2n-2j-1)!}{(j-i)!(j-1)!(n-j-i)!(n-j-1)!}\cdot p_i\\
& \equal{\eqref{E:binomial2}}
\begin{cases}
\sum_{i=0}^{\lfloor \frac{n}{2} \rfloor}  \frac{(n-2i)!}{n!} (-1)^i \frac{(2i-1)!!(n-1)!!}{(n-2i)!!}\cdot p_i & \text{ if $n$ is even}\\ 0 & \text{ if $n$ is odd} \end{cases}\\
& =\begin{cases}
\sum_{i=0}^{\lfloor \frac{n}{2} \rfloor} (-1)^i\frac{(2i-1)!!(n-2i-1)!!}{n!!}\cdot p_i
& \text{ if $n$ is even}\\ 0 & \text{ if $n$ is odd} \end{cases}\\
& \equal{\eqref{E:def-e}} 
 \begin{cases} \rho_{(n)}& \text{ if $n$ is even}\\ 0 & \text{ if $n$ is odd.} \end{cases}
\end{align*}
Equation~\eqref{E:e0torho0} can be deduced similarly from~\eqref{E:binomial1}.
\end{proof}

The dimensions of the left ideals of the group algebra $\field B_n$ generated by the idempotents $e_{(n)}$ and $\frac{1}{2}e_{(0,n)}$ are~\cite[Proposition 2.5 and page 108] {B92}
\[\dim (\field B_n) e_{(n)}= (2n-1)!! \text{ \ and \ } \dim (\field B_n) e_{(0,n)}=(2n-2)!!\,.\]
We calculate the dimensions of the left ideals of the group algebra $\field S_n$ generated by the idempotents $\rho_{(n)}$ and $\rho_{(0,n)}$.

\begin{proposition}\label{P:dim-ideals}
For even $n$,
\[\dim (\field S_n) \rho_{(n)}= (n-1)!!^2\,,\]
and for odd $n$,
\[\dim (\field S_n) \rho_{(0,n)}= (n-1)!\,.\]
\end{proposition}
\begin{proof} If $e$ is an idempotent of an algebra $A$ then $\dim Ae=\tr(r_e)$, the trace of the map $r_e:A\to A$, $r_e(a)=ae$ (since this is a projection onto $Ae$). If $A$ is a group algebra then $\tr(r_e)$ equals the coefficient of the identity of the group in $e$ times the order of the group. 
Assume $n$ is even, so $\rho_{(n)}$  is idempotent.
{}From~\eqref{E:def-rho} we see that the coefficient of the identity in $\rho_{(n)}$ is $\frac{(n-1)!!}{n!!}$, so 
\[\dim (\field S_n) \rho_{(n)}=\frac{(n-1)!!}{n!!}\cdot n!=(n-1)!!^2\,.\]
If $n$ is odd,  $\rho_{(0,n)}$  is idempotent, and by~\eqref{E:def-rho0} the coefficient of the identity in $\rho_{(0,n)}$ is $\frac{(n-2)!!}{n!!}$, so 
\[\dim (\field S_n) \rho_{(0,n)}=\frac{(n-2)!!}{n!!}\cdot n!=(n-2)!!(n-1)!!=(n-1)!\,.\]
\end{proof}

\subsection{A basis of semiidempotents}\label{S:semiidempotents}
 We build a basis of the peak algebra  by means of the convolution product.

\begin{definition}\label{D:semiidempotents}
For any  pseudocomposition $\beta=(b_0,b_1,\ldots,b_k)$ of $n$, let
\begin{equation}\label{E:ebasis}
e_\beta:=e_{(b_0)}\ast e_{(0,b_1)}\ast\cdots\ast e_{(0,b_k)}\,.
\end{equation}
Similarly, given an almost-odd composition $\gamma=(b_0,b_1,\ldots,b_k)$ of $n$, let
\begin{equation}\label{E:rhobasis}
\rho_\gamma:=\rho_{(b_0)}\ast \rho_{(0,b_1)}\ast\cdots\ast \rho_{(0,b_k)}\,.
\end{equation}
\end{definition}

Since $\Sol{B}$ is a graded right $\I^0$-module, $e_{\beta}\in\SolB$. 
Similarly, $\rho_{\gamma}\in\pppn$.

\begin{proposition}\label{P:etorho} Let $\beta$ be a pseudocomposition. Then
\begin{equation}\label{E:etorhotwo}
\varphi(e_{\beta}) =\begin{cases}
2^{k(\beta)}\rho_{\beta} & \text{ if $\beta$ is almost-odd,}\\ 
0 & \text{ otherwise.} \end{cases}
\end{equation}
\end{proposition}
\begin{proof} Since $\varphi$ preserves convolution products, 
\[\varphi(e_\beta)=\varphi(e_{(b_0)})\ast\varphi(e_{(0,b_1)})\ast\cdots\ast \varphi(e_{(0,b_k)})\,.\]
The result follows at once from Theorem~\ref{T:etorho}.
\end{proof}

The elements $e_\beta$ were introduced in~\cite{BB92b,B92} (where they are denoted $I_p$). It is shown in~\cite[page 106]{B92} that as $\beta$ runs over all pseudocompositions of $n$, the elements $e_{\beta}$ form a linear basis of $\SolB$. Moreover, each $e_\beta$ is a semiidempotent~\cite[Corollary 2.8]{B92}. 

\begin{corollary}\label{C:rhobasis} As $\gamma$ runs over the almost-odd compositions of $n$, the elements $\rho_\gamma$ form a basis of semiidempotents of $\pppn$.
\end{corollary}
\begin{proof} The surjectivity of $\varphi:\SolB\onto\pppn$ together with Proposition~\ref{P:etorho} imply that the elements $\rho_\gamma$ span $\pppn$. Since the dimension of $\pppn$ is the number of almost-odd compositions of $n$, they form a basis. Since each $e_\gamma$ is a semiidempotent, so is each $\rho_\gamma$.
\end{proof}

\subsection{Commutative semisimple subalgebras}\label{S:commutative}
Let $\sB$ denote the linear span of the elements $y_j$, $j=0,\ldots,n$, 
$i_n^0$ the linear span of the elements $y^0_j$, $j=1,\ldots,n$, and
\[\hsB:=\sB+i_n^0\,.\]
 It is known that $\hsB$ is a commutative semisimple subalgebra of $\SolB$ of dimension $2n$, $\sB$ is a subalgebra of $\hsB$ of
dimension $n+1$, and $i_n^0$ is an ideal of $\hsB$ of dimension $n$~\cite[Section 4.2]{Mah01}, ~\cite[Theorem 6.1]{ABN}.  Let
\begin{equation*}\label{E:def-x}
x_j  :=\sumsub{J\subseteq[0,n-1]\\\#J=j}X_J \text{ \ and \  }
x^0_j :=\sumsub{J\subseteq[0,n-1],\ 0\in J\\\#J=j}X_{J}\,.
\end{equation*}
The elements $x_j$, $j=0,\ldots,n$, form a basis of $\sB$, and the elements
$x^0_j$, $j=1,\ldots,n$, form a basis of $i_n^0$~\cite[Section 6.1]{ABN}. The idempotents $e_{(n)}$ and $e_{(0,n)}$ can be expressed in these bases as follows:
\[e_{(n)}=\sum_{j=0}^{n} (-1)^j\frac{(2j-1)!!}{(2j)!!}x_j
\text{ \ and \ }e_{(0,n)}=\sum_{j=1}^{n} (-1)^{j-1}\frac{1}{j}x^0_j\,.\]
These formulas can be found in~\cite[Section 2]{B92} or~\cite[Section 4.2]{Mah01}.

We discuss  peak analogs of these formulas.
Let  $\wp_n$ denote the linear span of the elements $p_j$, $j=0,\ldots,\ipartn$,
$\wp^0_n$ the linear span of the elements  $p_j^0$, $j=1,\ldots,\ipart{n+1}$,
and 
\[\hwp_n:=\wp_n+\wp^0_n\,.\]
 We know  that $\hwp_n$ 
 is a commutative semisimple subalgebra of $\pppn$ of dimension $n$, $\wp_n$ is a subalgebra of dimension $\ipartn+1$, and $\wp^0_n$ is an ideal of dimension $\ipart{n+1}$~\cite[Theorem 6.8]{ABN}. Define elements
\begin{equation}\label{E:def-q}
q_j:=\sumsub{F \text{ sparse}\\ \#F=j}Q_F \text{ \ and \ }
q^0_j:=\sumsub{F \text{ sparse, }1\notin F\\ \#F=j-1}Q_F\,.
\end{equation}

\begin{proposition}\label{P:qinp}
\begin{equation}\label{E:qinp}
q_j=\sum_{i=j}^{\ipartn}\binom{i}{j}p_i \text{ \ and \ }
q_j^0=\sum_{i=j}^{\ipart{n+1}}\binom{i-1}{j-1}p_i^0\,.
\end{equation}
\end{proposition}
\begin{proof} We have
\begin{align*}
q_j & =\sumsub{F\subseteq G \text{ sparse}\\ \#F=j}P_G=
\sum_{G \text{ sparse}}\#\{F\subseteq G \mid\ \#F=j\}P_G\\
\notag&=\sum_{G \text{ sparse}}\binom{\#G}{j}P_G=
\sum_{i=j}^{\ipartn}\binom{i}{j}p_i\,.
\end{align*}
The formula for $q_j^0$ is similar.
\end{proof}

It follows that the $q_j$, $j=0,\ldots,\ipartn$, form a basis of $\wp_n$ and
 the $q^0_j$, $j=1,\ldots,\ipart{n+1}$, form a basis of $\wp^0_n$. The elements $\rho_{(n)}$ and $\rho_{(0,n)}$ can be expressed in these bases as follows.
\begin{proposition}\label{P:rhoinq} For every $n$,
\[\rho_{(n)}=\sum_{i=0}^{\ipartn} (-1)^i\frac{(n-2i-1)!!}{(n-2i)!!}q_i
\text{ \ and \ }\rho_{(0,n)}=\sum_{i=1}^{\ipart{n+1}} (-1)^{i-1}\frac{1}{n-2i+2}q^0_i\,.\]
\end{proposition}
\begin{proof} Left to the reader. 
\end{proof}

The canonical map $\varphi:\SolB\to\pppn$ admits the following expressions on the bases $x_j$ and $q_j$.
 \begin{proposition}\label{P:xtoq}
\[\varphi(x_j)=2^j\sum_{i=0}^{\ipartn} (-1)^i\binom{n-2i}{j}q_i 
\text{ \ and \ }\varphi(x^0_j)=2^j\sum_{i=1}^{\ipart{n+1}} (-1)^{i-1}\binom{n+1-2i}{j-1}q^0_i \,.\]
\end{proposition}
\begin{proof} We have
\[\varphi(x_j)=\sumsub{J\subseteq[0,n-1]\\\#J=j}\varphi(X_J)\ \ \equal{\eqref{E:XinQ}}
\ \ 2^{j}\cdot\!\!\!\!\!\!\!\!\!\!\!\!\sumsub{G \text{ sparse}\\ G\subseteq [n-1]\setminus\bigl(J\cup(J+1)\bigr)} \!\!\! (-1)^{\#G}Q_G\,.\]
As seen in the proof of Lemma~\ref{L:sparse-almostodd},
 $G\subseteq [n-1]\setminus\bigl(J\cup(J+1)\bigr)\iff
J\subseteq [0,n-1]\setminus\bigl(G\cup(G-1)\bigr)$. 
Once a sparse subset $G$ has been fixed, there are $\binom{n-2\#G}{j}$  subsets $J$ satisfying this condition, since $G$ and $G-1$ are disjoint. The formula for $\varphi(x_j)$ follows. The formula for $\varphi(x^0_j)$ can be derived similarly.
\end{proof}

Consider the maps $\beta:\SolB\to\Sol{B_{n-2}}$ and $\pi:\pppn\to\ppp{n-2}$ of
Section~\ref{S:chains}. We have~\cite[Proposition 6.10]{ABN}
\begin{align*}
\beta(x_j) & = \begin{cases} x_j & \text{ if } 0\leq j<n\,,\\
                          0 & \text{ if } j=n\,.
             \end{cases}\\
\intertext{Similarly, one sees that}
\pi(q_j) & =\begin{cases} -q_{j-1} & \text{ if $j=1,\ldots,\ipartn$,}\\
0 & \text{ if $j=0$.}\end{cases} 
\end{align*}
Since $e_{(0,n)}\in\I_n^0=\ker(\beta)$, we have $\beta(e_{(0,n)})=0$ for every $n$. Similarly,  $\pi(\rho_{(0,n)})=0$ for every $n$.
 On the other hand,
\begin{proposition}\label{P:beta-e}
For every $n$,
\[\beta(e_{(n)})=e_{(n-1)}
\text{ \ and \ }
\pi(\rho_{(n)})=\rho_{(n-2)}\,.\]
\end{proposition}
\begin{proof}
These follow easily from the preceding formulas.
\end{proof}

\section{The action on Lie monomials}\label{S:Lie}

\subsection{Preliminaries}\label{S:preliminaries}
Let $T(V)=\oplus_{n\geq 0}V^{\otimes n}$ be the tensor algebra of a vector space $V$. It is a Hopf
algebra with coproduct determined by $\Delta(v)=1\otimes v+v\otimes 1$ for all $v\in V$.
A {\em Lie polynomial} is a primitive element of $T(V)$.
A {\em Lie monomial} is a product of Lie polynomials. View $T(V)$ as a Lie algebra under the commutator bracket $[a,b]=ab-ba$. The subspace $L(V)$
of Lie polynomials may also be described as the Lie subalgebra of $T(V)$
generated by $V$. This turns out to be the free Lie algebra generated by $V$.
We have $L(V)=\oplus_{n\geq 1}L_n(V)$ with $L_n(V)=L(V)\cap V^{\otimes n}$.

Recall the right action of $S_n$ on the tensor power $V^{\otimes n}$ of a  vector space $V$:
\[(v_1\ldots v_n)\cdot\sigma = v_{\sigma(1)}\ldots v_{\sigma(n)}\,.\]
There is also a left action of $GL(V)$ on $V^{\otimes n}$ given by
\[g\cdot(v_1\ldots v_n)=(g\cdot v_1)\ldots(g\cdot v_n)\,.\]
These actions commute: for any $g\in GL(V)$, $a\in V^{\otimes n}$, and $\sigma\in S_n$,
\[g\cdot (a\cdot\sigma)=(g\cdot a)\cdot \sigma\,.\]
A classical result (Schur-Weyl duality) states that  if $\dim V\geq n$ then $\field S_n$ may be recovered as those endomorphisms of $V^{\otimes n}$ which
commute with the action of $GL(V)$.

Similarly, an important result of Garsia and Reutenauer characterizes
which elements of the group algebra $\field S_n$ belong to the descent algebra $\SolA$ in terms of their action on  Lie monomials
~\cite[Theorem 4.5]{GReu}. Their result may be stated as follows. An element $\phi\in\field S$ belongs to $\Sol{A}$ if and only if for every Lie polynomials
$p_1,\ldots,p_k\in L(V)$, the subspace
\begin{equation}\label{E:action-GR}
\Span\{p_{s(1)}\ldots p_{s(k)}  \mid s\in S_k\}\subseteq T(V)
\end{equation}
is invariant under the action of $\phi$.

Schocker obtained a 
characterization for the elements of the peak ideal $\pppn^0$
in terms of the action on Lie monomials~\cite[Main Theorem 8]{Sch}. 
Let $L(V)=L_e(V)\oplus L_o(V)$ denote the decomposition into polynomials
of even and odd degrees, i.e., $L_e(V)=\oplus_{n \text{ even}}L_n(V)$, and  $L_o(V)=\oplus_{n \text{ odd}}L_n(V)$.
Schocker's result states that  an element $\phi\in\Sol{A}$ belongs to $\ppp{}^0$ if and only if for every Lie polynomials $p_1,\ldots,p_k\in L(V)$ with $p_1\in L_e(V)$,
\begin{equation}\label{E:action-Schocker}
(p_1\ldots p_k)\cdot\phi=0\,.
\end{equation}

Below we present a characterization for the elements of the
peak algebra $\pppn$ that is analogous to that of Garsia and Reutenauer, both in content and proof (Theorem~\ref{T:charpeak}).  Furthermore, we provide
a characterization for the elements of
each ideal $\pppn^j$ that interpolates between Schocker's
characterization of the peak ideal and our characterization
of the peak algebra (Theorem~\ref{T:charpeakideal}).

\subsection{The action of signed permutations}\label{S:action-signed}

Suppose the vector space $V$ is endowed with an involution
\[v\mapsto \inv{v}\,,\quad \inv{\inv{v}}=v\,.\]
We extend the involution to  $T(V)$ by
\[\inv{v_1\ldots v_n}:=\inv{v_n}\ldots\inv{v_1}\,.\]
Thus $a\mapsto\inv{a}$ is an anti-automorphism of algebras of $T(V)$.
We say that an element $a\in T(V)$ is {\em invariant} if
$\inv{a}=a$, and {\em skew-invariant} if $\inv{a}=-a$. We obtain decompositions
\[T(V)=T_i(V)\oplus T_s(V) \text{ \ and \ }L(V)=L_i(V)\oplus L_s(V)\]
into invariants and skew-invariants elements.

The group $B_n$ acts on $V^{\otimes n}$ via
\begin{equation}\label{E:Baction}
(v_1\ldots v_n)\cdot \sigma=v_{\sigma(1)}^{\pm}\ldots v_{\sigma(n)}^{\pm}\,,
\text{ \ where \ }
v_{\sigma(i)}^{\pm}:=\begin{cases} v_{\sigma(i)} & \text{ if }\sigma(i)>0\,,\\
\inv{v_{-\sigma(i)}} & \text{ if }\sigma(i)<0\,.\end{cases}
\end{equation}

Let $\iota_n:\field B_n\to\End(V^{\otimes n})$ be $\iota(\sigma)(a)=a\cdot\sigma$. Summing over $n$ we get a map
\[ \iota: \field B\to\End(T(V))\,.\]
The external product of $\field B$ (Section~\ref{S:convolution}) corresponds to the
convolution of endomorphisms under $\iota$: for any $\sigma\in B_p$ and $\tau\in B_q$,
\begin{equation}\label{E:convolutiontwo}
\iota(\sigma\ast\tau)=m\bigl(\iota(\sigma)\otimes\iota(\tau)\bigr)\Delta
\end{equation}
where $m$ and $\Delta$ are the product and coproduct of $T(V)$.

Consider the operator $\nabla:T(V)\to T(V)$ defined by
\[\nabla(a)= a+\inv{a}\,.\]
The following result is central for our purposes. It generalizes~\cite[Proposition 8.8]{ABN}.

\begin{proposition}\label{P:Xaction} Let $p_1,\ldots,p_k$ be homogeneous Lie polynomials and $n=\sum_{i=1}^k \deg(p_i)$. Then
\begin{equation}\label{E:Xaction}
(p_1\ldots p_k)\cdot X_{(0,n)}=\nabla\Bigl(\ldots\nabla\bigl(\nabla(p_1)p_2\bigr)
p_3\ldots p_k\Bigr)\,.
\end{equation}
In particular, if $p_1$ is skew-invariant, then
\begin{equation}\label{E:Xactionskew}
(p_1\ldots p_k)\cdot X_{(0,n)}=0\,.
\end{equation}
\end{proposition}
\begin{proof} For $p\geq 0$, let $1_p:=12\ldots p\in B_p$ denote the identity permutation and let
$\bar{1}_p:=\bar{p}\ldots\bar{2}\bar{1}\in B_p$. Note that $\bar{1}_p(a)=\inv{a}$ for any $a\in V^{\otimes p}$. Let $R_{(p,q)}:=\bar{1}_p\ast 1_q\in\field B_{p+q}$.  We have
\[X_{(0,n)}=\sum_{p=0}^nR_{(p,n-p)}\]
(see the proof of Proposition 7.13 in ~\cite{ABN} for a more general result).
We make use of~\eqref{E:convolutiontwo} to analyze the action of $X_{(0,n)}$.
Since each $p_i$ is primitive, we have
\[\Delta(p_1\ldots p_k)=\sum_{S\sqcup T=[k]}p_S\otimes p_T\,,\]
where, if $S=\{i_1<\ldots<i_h\}$, then $p_S:=p_{i_1}\ldots p_{i_k}$.
Therefore,
\[(p_1\ldots p_k)\cdot X_{(0,n)}=\sum_{S\sqcup T=[k]}\inv{p_S}p_T\,.\]
 We verify that
\[\sum_{S\sqcup T=[k]}\inv{p_S}p_T=
\nabla\Bigl(\ldots\nabla\bigl(\nabla(p_1)p_2\bigr)p_3\ldots p_k\Bigr)\]
 by induction on $k$. If $k=1$, both sides equal $p_1+\inv{p_1}$. Assume the result holds for $k-1$. Then the right hand side equals
 \[\nabla\Bigl(\sum_{S\sqcup T=[k-1]}\inv{p_S}p_Tp_k\Bigr)=\sum_{S\sqcup T=[k-1]}\inv{p_S}p_Tp_k+\sum_{S\sqcup T=[k-1]}\inv{p_k}\inv{p_T}p_S\,.\]
 The subsets $S$ from the first sum, together with the subsets $T\cup\{k\}$ from the second sum, traverse all subsets of $[k]$, and we obtain  the left hand side.
\end{proof}

\subsection{The action of elements of the peak algebra}\label{action-peak}

We revert to the case of an arbitray vector space $V$. We endow it with the
trivial involution $\inv{v}:=v$. The induced involution on $T(V)$ is
\[\inv{v_1\ldots v_n}= v_n\ldots v_1\,.\]
There are invariants and skew-invariants of arbitrary degrees. However, a Lie polynomial is invariant (skew-invariant) if and only if it is odd (even).
\begin{lemma}\label{L:even-odd} We have
\[L_i(V)=L_o(V) \text{ \ and \ }L_s(V)=L_e(V)\,.\]
\end{lemma}
\begin{proof} We have $L(V)=L_i(V)\oplus L_s(V)=L_o(V)\oplus L_e(V)$,
so it suffices to show that $L_o(V)\subseteq L_i(V)$ and $L_e(V)\subseteq L_s(V)$. We verify the first inclusion, the second is similar. Let $p\in L_o(V)$.
We may assume that $p$ is homogeneous and we argue by induction on its
degree. If $\deg(p)=1$ we have $\inv{p}=p$ because the involution is trivial on $V$. If $\deg(p)>1$ then $p$ is a linear combination of polynomials of the form $[a,b]$, with $a$ and $b$ homogeneous Lie polynomials of smaller degree.
Since $\deg(p)$ is odd, one of the polynomials $a$ and $b$ is even and the other is odd.  By induction hypothesis, one of them is skew-invariant and the other is invariant. Hence,
\[\inv{[a,b]}=\inv{ab-ba}=\inv{b}\inv{a}-\inv{a}\inv{b}=-ba+ab=[a,b]\,.\]
Thus $[a,b]$, and hence $p$, is invariant.
\end{proof}

Observe that, since the involution is trivial on $V$, the action of $B_n$ on $V^{\otimes n}$~\eqref{E:Baction} is
\[(v_1\ldots v_n)\cdot \sigma=v_{\abs{\sigma(1)}}\ldots v_{\abs{\sigma(n)}}\,.\]
Therefore, the action of $B_n$ descends to the (usual) action of $S_n$ on $V^{\otimes n}$ via the canonical map $\varphi:B_n\to S_n$. 

Using results of Bergeron and Bergeron, we may now describe  the action on the tensor algebra of the idempotents $\rho_{(n)}\in\pppn$ and $\rho_{(0,n)}\in\pppn^0$  of Theorem~\ref{T:etorho}. The latter acts  as the projection onto the subspace of odd Lie polynomials, the former as the 
projection onto the subalgebra generated by even Lie polynomials.
\begin{lemma}\label{L:rho-lie} 
\begin{align*}
T(V)\cdot\sum_{n \text{ even}} \rho_{(n)} & = \text{ the subalgebra of $T(V)$ generated by $L_e(V)$,}\\
T(V)\cdot\sum_{n \text{ odd}} \rho_{(0,n)} & = L_o(V)\,.
\end{align*}
\end{lemma}
\begin{proof} According to~\cite[Theorem 2.2]{B92}, the element $\sum_n e_{(n)}$ projects $T(V)$ onto the subalgebra of $T(V)$ generated by $L_s(V)$. According to~\cite[Theorem 2]{BB92a} or~\cite[Theorem 2.1]{B92},
 the element $\sum_n e_{(0,n)}$ projects $T(V)$ onto $L_i(V)$. Together with Theorem~\ref{T:etorho} and Lemma~\ref{L:even-odd} this gives the result.
\end{proof}

The sum of all permutations in $S_n$ with no interior peaks is
\[P_{(0,n)}:=P_\emptyset+P_{\{1\}}\in\pppn\,.\]
According to ~\eqref{E:XinP} and~\eqref{E:XontoTO},
\begin{equation}\label{E:Pempty}
P_{(0,n)}=\frac{1}{2}\varphi(X_{(0,n)}) \text{ for any $n$, and }
P_{(0,n)}=\TO_{(0,n)} \text{ if $n$ is odd.} 
\end{equation}

Consider  the operator $T(V)\times T(V)\to T(V)$ defined on homogeneous elements $a$ and $b$ by
\[ \{a,b\}=ab+(-1)^{\deg(b)-1}ba\,.\]

The following result  describes the action of $P_{(0,n)}$ on Lie monomials. It generalizes~\cite[Proposition 8.9]{ABN} and is closely related to~\cite[Lemma 5.11]{KLT}.
\begin{proposition}\label{P:Xaction-peaks} Let $p_1,\ldots,p_k$ be homogeneous Lie polynomials and $n=\sum_{i=1}^k\deg(p_i)$. Then
\begin{equation}\label{E:Xaction-peaks}
(p_1\ldots p_k)\cdot P_{(0,n)}=\begin{cases}
\Bigl\{\ldots\bigl\{\{p_1,p_2\},p_3\bigr\},\ldots,p_k\Bigr\} & \text{ if $p_1$ is odd,}\\
0 & \text{ if $p_1$ is even,} \end{cases}
\end{equation}
\end{proposition}
\begin{proof} {}From~\eqref{E:Xaction} and~\eqref{E:Pempty} we get
\[(p_1\ldots p_k)\cdot P_{(0,n)}=
\frac{1}{2}\nabla\Bigl(\ldots\nabla\bigl(\nabla(p_1)p_2\bigr)p_3\ldots p_k\Bigr)\,.\]
If $p_1$ is even then $\nabla(p_1)=p_1+\inv{p_1}=0$ by 
Lemma~\ref{L:even-odd}, and we are done.

When $p_1$ is odd we argue by induction on $k$.  Let $\eta_k=\nabla\Bigl(\ldots\nabla\bigl(\nabla(p_1)p_2\bigr)p_3\ldots p_k\Bigr)$ and
$\theta_k=\Bigl\{\ldots\bigl\{\{p_1,p_2\},p_3\bigr\},\ldots,p_k\Bigr\}$. We have to show that $\eta_k=2\theta_k$.

For $k=1$ we have, by Lemma~\ref{L:even-odd},
$\eta_1=\nabla(p_1)=p_1+\inv{p_1}=2p_1=2\theta_1$,
so the result holds. 

Assume the result holds for $k$. Note that for any homogeneous Lie polynomial $p$ we have $\inv{p}=(-1)^{\deg(p)-1}p$, by Lemma~\ref{L:even-odd}. Also,
since $\eta_k$ is in the image of $\nabla$, $\inv{\eta_k}=\eta_k$. Hence,
\begin{align*}
\eta_{k+1} &=\nabla(\eta_k p_{k+1})=\eta_k p_{k+1}+\inv{p_{k+1}}\,\inv{\eta_k}=
\eta_k p_{k+1}+ (-1)^{\deg(p_{k+1})-1}p_{k+1}\eta_k\\
&=2\theta_k p_{k+1}+ (-1)^{\deg(p_{k+1})-1}p_{k+1}2\theta_k=2\{\theta_k, p_{k+1}\}=2\theta_{k+1}\,,
\end{align*}
as needed.
\end{proof}

We may now derive the characterization of the peak algebra in terms of the
action on Lie monomials.

\begin{theorem}\label{T:charpeak} Let $V$ be an infinite-dimensional vector space. An element $\phi\in\field S$ belongs to $\ppp{}$ if and only if for every Lie polynomials
$p_1,\ldots,p_u\in L_e(V)$ and $q_1,\ldots,q_v\in L_o(V)$, the subspace
\begin{equation}\label{E:subspace}
\Span\{p_1\ldots p_u q_{s(1)}\ldots q_{s(v)}  \mid s\in S_v\}\subseteq T(V)
\end{equation}
is invariant under the action of $\phi$.
\end{theorem}
\begin{proof} We first show that  subspace~\eqref{E:subspace} is invariant under any element $\phi$ of the peak algebra. We may assume that $p_i$, $q_j$ are homogeneous Lie polynomials and
$\phi=\TO_\beta$, $\beta=(b_0,b_1,\ldots,b_k)$ an almost-odd composition of $n:=\sum_{i=1}^u \deg(p_i)+\sum_{j=1}^v \deg(q_j)$. We have
$\TO_{(b_0,b_1,\ldots,b_k)} =\TO_{(b_0)}\ast \TO_{(0,b_1)}\ast\cdots\ast \TO_{(0,b_k)}$~\eqref{E:TOfree}. For $i=1,\ldots,u+v$, let
\[\ell_i=\begin{cases} p_i & \text{ if $1\leq i\leq u$,}\\
q_{i-u} & \text{ if $u+1\leq i\leq u+v$.} \end{cases}\]
Since Lie polynomials are primitive elements,
\[\Delta^{(k)}(p_1\ldots p_u q_{1}\ldots q_{v})=\sum_{T_0\sqcup\cdots\sqcup T_k=[u+v]}\ell_{T_0}\otimes\cdots\otimes\ell_{T_k}\,,\]
where $\ell_T:=\prod_{i\in T}\ell_i$ (product in {\em increasing} order of the indices, as in the proof of Proposition~\ref{P:Xaction}). Therefore, by~\eqref{E:convolutiontwo},
\[(p_1\ldots p_u q_{1}\ldots q_{v})\cdot\TO_\beta=\sumsub{T_0\sqcup\cdots\sqcup T_k=[u+v]\\\deg(\ell_{T_i})=b_i}
\bigl(\ell_{T_0}\cdot \TO_{(b_0)}\bigr)\bigl(\ell_{T_1}\cdot\TO_{(0,b_1)}\bigr)\ldots\bigl(\ell_{T_k}\cdot\TO_{(0,b_k)}\bigr)\,.\]
In this sum, if for any $i\geq 1$ the subset $T_i$ contains an element from $[u]$,
then the first factor of the Lie monomial $\ell_{T_i}$ is an even Lie polynomial
(one of the $p$'s); then, by~\eqref{E:Pempty} and~\eqref{E:Xaction-peaks},
\[\ell_{T_i}\cdot\TO_{(0,b_i)}=0\,.\]
Thus the only terms that contribute to this sum are those for which $T_0\supseteq [u]$. In this case, since the element $\TO_{(b_0)}$ is the identity of $S_{b_0}$, we have 
\[\ell_{T_0}\cdot \TO_{(b_0)}=\ell_{T_0}=p_1\ldots p_uq_1\ldots q_{\#T_0-u}\,.\]
On the other hand, the elements $\TO_{(0,b_i)}$ are, in particular, elements of the descent algebra $\SolA$, so by the result of Garsia and Reutenauer~\eqref{E:action-GR} each 
$\ell_{T_i}\cdot\TO_{(0,b_i)}$ is a linear combination of Lie monomials of the form $\ell_{s(j_1)}\ldots\ell_{s(j_{v_i})}$, as $s$ runs over the permutations of the set $T_i:=\{j_1,\ldots,j_{v_i}\}$. It follows that $(p_1\ldots p_u q_{1}\ldots q_{v})\cdot\TO_\beta$ is a linear combination of Lie monomials of the form
\[p_1\ldots p_uq_{s(1)}\ldots q_{s(v)}\,,\]
as $s$ runs over the permutations of $[v]$. This proves the invariance of 
subspace~\eqref{E:subspace}.

We now prove the converse. Start from an element $\phi\in\field S_n$ under 
whose action any subspace~\eqref{E:subspace} is invariant. 

Fix an almost-odd composition $\beta=(b_0,b_1,\ldots,b_k)$ of $n$.
Let $I_0\sqcup I_1\sqcup\cdots\sqcup I_k=[n]$ be the decomposition of
$[n]$ into consecutive segments of lengths $b_0,b_1,\ldots,b_k$. Thus
$I_0=\{1,\ldots,b_0\}$, $I_1=\{b_0+1,\ldots,b_0+b_1\}$, etc.

Let $v_1,\ldots,v_n$ be linearly independent elements of $V$.  Define
\[P:=v_{I_0}\cdot \rho_{(b_0)}\,,\ 
q_1 :=v_{I_1}\cdot \rho_{(0,b_1)}\,,\ \ldots,\ 
q_k :=v_{I_k}\cdot \rho_{(0,b_k)}\,,\]
where $\rho_{(n)}$ and $\rho_{(0,n)}$ are the idempotents of Theorem~\eqref{T:etorho}, and in each $v_T:=\prod_{i\in T}v_i$  the product is in  increasing order of the indices (as before). 

 By Lemma~\ref{L:rho-lie}, $q_1,\ldots,q_k\in L_o(V)$, and $P$ belongs to the subalgebra generated by $L_e(V)$, so there are even Lie polynomials $p_1,\ldots,p_h\in L_e(V)$ such that $P:=p_1\ldots p_h$.
Therefore, our hypothesis implies the existence of scalars $c_s$ indexed by permutations $s\in S_k$ such that
\begin{equation}\tag{$*$}
(Pq_1\ldots q_k)\cdot\phi=\sum_{s\in S_k} c_s Pq_{s(1)}\ldots q_{s(k)}\,.
\end{equation}
 Fix a decomposition $T_0\sqcup T_1\sqcup\cdots\sqcup T_k=[n]$ with $\#T_i=b_i$ for every $i$. Let $\gamma$ be the unique permutation of $[n]$
 such that $\gamma(I_i)=T_i$ and $\gamma$ restricted to each $I_i$ is order-preserving. 
 Since the $v_i$ are linearly independent, there is  a linear transformation $g\in GL(V)$ such that $g(v_i)=v_{\gamma(i)}$ for each $i$. Note that $g\cdot v_{I_i}=v_{T_i}$. Since the actions of $GL(V)$ and $S_n$ on $V^{\otimes n}$ commute with each other, we have $g\cdot P=(g\cdot v_{I_0})\cdot\rho_{(b_0)}=v_{T_0}\cdot\rho_{(b_0)}$, and similarly
 $g\cdot q_i=v_{T_i}\cdot\rho_{(0,b_i)}$. Thus, 
  acting with $g$ from the left  on both sides of $(*)$ we obtain 
 \begin{multline*}
\bigl((v_{T_0}\cdot \rho_{(b_0)})(v_{T_1}\cdot \rho_{(0,b_1)})\ldots (v_{T_k}\cdot \rho_{(0,b_k)})\bigr)\cdot\phi\\
=\sum_{s\in S_k} c_s (v_{T_{0}}\cdot \rho_{(b_0)})(v_{T_{s(1)}}\cdot \rho_{(0,b_{s(1)})})\ldots (v_{T_{s(k)}}\cdot \rho_{(0,b_{s(k)})})\,.
\end{multline*}
Note that the coefficients $c_s$ are the same for all decompositions $\{T_i\}$.
Summing over all such decompositions, we obtain
 \begin{multline*}
(v_1\ldots v_n)\cdot(\rho_{(b_0)}\ast\rho_{(0,b_1)}\ast\cdots\ast\rho_{(0,b_k)})\cdot\phi\\
=\sum_{s\in S_k} c_s (v_1\ldots v_n)\cdot (\rho_{(b_0)}\ast\rho_{(0,b_{s(1)})}\ast \cdots\ast \rho_{(0,b_{s(k)})})\,.
\end{multline*}
 (The convolution product gives rise to the sum over those decompositions because the $v_i$ are primitive elements.) Now, by~\eqref{E:rhobasis}, this equation may be rewritten as
\[(v_1\ldots v_n)\cdot(\rho_{(b_0,b_1,\ldots,b_k)})\cdot\phi
=\sum_{s\in S_k} c_s (v_1\ldots v_n)\cdot \rho_{(b_0,b_{s(1)},\ldots,b_{s(k)})}\,.\]
Since the $v_i$ are linearly independent, this implies
\[\rho_{(b_0,b_1,\ldots,b_k)}\cdot\phi
=\sum_{s\in S_k} c_s  \rho_{(b_0,b_{s(1)},\ldots,b_{s(k)})}\,.\]
In particular, for any almost-odd composition $\beta$,
\[\rho_{\beta}\cdot\phi\in\pppn\,.\]
Since the $\rho_\beta$ form a basis of $\pppn$ (Corollary~\ref{C:rhobasis}),
we may write $1\in\pppn$ as a linear combination of these elements, and conclude that $\phi\in\pppn$. This completes the proof.
\end{proof}

\begin{example}\label{Ex:charpeak}
Let $a,b,c,d\in V$ and consider the Lie polynomials $p_1=[a,b]$, $p_2=c$, and $p_3=[a,[b,d]]$. We have
\[p_1p_2p_3=abcabd-abcadb-abcbda+abcdba-bacabd+bacadb+bacbda-bacdba\,.\]
The total degree is $n=6$. The action of
\[P_{\{5\}}=123465+123564+124563+134562+234561\in\ppp{6}\]
may be explicitly calculated as follows:
\begin{align*}
(p_1p_2p_3)\cdot 123465&= abcadb-abcabd-abcbad+abcdab-bacadb+bacabd+bacbad-bacdab\\
(p_1p_2p_3)\cdot 123564&= abcbda-abcdba-abcdab+abcbad-bacbda+bacdba+bacdab-bacbad\\
(p_1p_2p_3)\cdot 124563&= ababdc-abadbc-abbdac+abdbac-baabdc+baadbc+babdac-badbac\\
(p_1p_2p_3)\cdot 134562&= acabdb-acadbb-acbdab+acdbab-bcabda+bcadba+bcbdaa-bcdbaa\\
(p_1p_2p_3)\cdot 234561&= bcabda-bcadba-bcbdaa+bcdbaa-acabdb+acadbb+acbdab-acdbab\,.
\end{align*}
It follows that
\begin{align*}
(p_1p_2p_3)\cdot P_{\{5\}}&= abcadb-abcabd 
\phantom{-abcbad+abcdab\,}-bacadb+bacabd
\phantom{+bacbad-bacdab}\\
&+ abcbda-abcdba\phantom{-abcdab+abcbad\ }-bacbda+bacdba\phantom{+bacdab-bacbad}\\
&+ ababdc-abadbc-abbdac+abdbac-baabdc+baadbc+babdac-badbac\\
&=-p_1p_2p_3+p_1p_3p_2\,.
\end{align*}
\end{example}

Theorem~\ref{T:charpeak} still holds if we only assume $\dim V\geq n$, provided we start from an element $\phi\in\field S_n$ (with the same proof).
Clearly one need only consider Lie monomials of total degree $n$ in~\eqref{E:subspace}.

\medskip

To derive the characterization of the ideals $\pppn^j$ in terms of the action
on Lie monomials, we analyze the behavior of the map $\pi:\pppn\to\ppp{n-2}$~\eqref{E:defpi} with respect to this action.

\begin{lemma}\label{L:action-pi} Let $\ell_0$ a Lie polynomial of degree $2$,  $\ell_1,\ldots,\ell_v$ homogeneous Lie polynomials, and
$n=2+\sum_{i=1}^v\deg(\ell_i)$. Then, for any $\phi\in\pppn$,
\begin{equation}\label{E:action-pi}
(\ell_o\ell_1\ldots \ell_v)\cdot\phi=\ell_0\bigl((\ell_1\ldots \ell_v)\cdot\pi(\phi)\bigr)\,.
\end{equation}
\end{lemma}
\begin{proof} It suffices to consider the case when $\phi=\TO_\beta$,
 $\beta=(b_0,b_1,\ldots,b_k)$ an almost-odd composition of $n$.
As in the proof of Theorem~\ref{T:charpeak},
\[(\ell_0 \ell_{1}\ldots \ell_{v})\cdot\TO_\beta=\sumsub{T_0\sqcup\cdots\sqcup T_k=[0,v]\\0\in T_0,\,\deg(\ell_{T_i})=b_i}
\bigl(\ell_{T_0}\cdot \TO_{(b_0)}\bigr)\bigl(\ell_{T_1}\cdot\TO_{(0,b_1)}\bigr)\ldots\bigl(\ell_{T_k}\cdot\TO_{(0,b_k)}\bigr)\,.\]
Note that if $0\in T_0$ then the first factor in the Lie monomial $\ell_{T_0}$  is $\ell_0$, and $\deg(\ell_{T_0})\geq 2$. Therefore,
If $b_0=0$, then no decomposition satisfies both
$0\in T_0$ and $\deg(\ell_{T_0})=b_0$, so $(\ell_0 \ell_{1}\ldots \ell_{v})\cdot\TO_\beta=0$. This agrees with the right hand side of~\eqref{E:action-pi}, because in this case  $\pi(\TO_\beta)=0$ by~\eqref{E:pi-TO}.

Assume $b_0\geq 2$. Since $\TO_{(b_0)}$ is the identity of $S_{b_0}$, we may write
\[(\ell_0 \ell_{1}\ldots \ell_{v})\cdot\TO_\beta=\sumsub{T_0\sqcup\cdots\sqcup T_k=[v]\\\deg(\ell_{T_0})=b_0-2,\,\deg(\ell_{T_i})=b_i}
\ell_0\bigl(\ell_{T_0}\cdot \TO_{(b_0-2)}\bigr)\bigl(\ell_{T_1}\cdot\TO_{(0,b_1)}\bigr)\ldots\bigl(\ell_{T_k}\cdot\TO_{(0,b_k)}\bigr)\,.\]
The right hand side equals
\[\ell_0\bigl((\ell_1\ldots \ell_v)\cdot\TO_{(b_0-2,b_1,\ldots,b_k)}\bigr)=
\ell_0\bigl((\ell_1\ldots \ell_v)\cdot\pi(\TO_\beta)\bigr)\]
by~\eqref{E:pi-TO}.
\end{proof}

\begin{theorem}\label{T:charpeakideal} Let $V$ be a vector space of dimension $\geq n$. Let $j=0,\ldots,\ipartn$. An element $\phi\in\field S_n$ belongs to $\pppn^j$
 if and only if for any homogeneous  Lie polynomials $p_1,\ldots,p_u\in L_e(V)$ and $q_1,\ldots,q_v\in L_o(V)$ we have that
\begin{equation}\label{E:action-us}
(p_1\ldots p_u q_1\ldots q_v)\cdot\phi = 
\begin{cases} 
{\displaystyle \sum_{s\in S_v}}c_s 
p_1\ldots p_u q_{s(1)}\ldots q_{s(v)} & \text{ if }
2j\geq \sum_{i=1}^u\deg(p_i)\,,\\
0 & \text{ if }2j<\sum_{i=1}^u\deg(p_i)\,,
\end{cases}
\end{equation}
where $c_{s}\in\field$ can be arbitrary scalars.
\end{theorem}
\begin{proof} Fix $j$ and suppose $\phi\in\pppn^j$.  Choose Lie polynomials as in~\eqref{E:action-us}. We may assume $\sum_{i=1}^u\deg(p_i)+\sum_{i=1}^v\deg(q_i)=n$ and $\phi=\TO_\beta$,
with $\beta=(b_0,b_1,\ldots,b_k)$ an almost-odd composition of $n$ with
$b_0\leq 2j$ (Corollary~\ref{C:peakideals}).
As in the proof of Theorem~\ref{T:charpeak}, we have
\begin{equation}\tag{$*$}
(p_1\ldots p_u q_{1}\ldots q_{v})\cdot\TO_\beta=\sumsub{T_0\sqcup\cdots\sqcup T_k=[u+v]\\ \deg(\ell_{T_i})=b_i} 
\bigl(\ell_{T_0}\cdot\TO_{(b_0)}\bigr)
\bigl(\ell_{T_1}\cdot\TO_{(0,b_1)}\bigr)\ldots\bigl(\ell_{T_k}\cdot\TO_{(0,b_k)}\bigr)\,,
\end{equation}
and the only decompositions $\{T_i\}$ that contribute to this sum have $T_0\supseteq[u]$. This condition implies that $p_{[u]}$ is a factor of
$\ell_{T_0}$ and $\deg(\ell_{T_0})\geq \deg(p_{[u]})$. 
If $2j<\deg(p_{[u]})$, then $b_0<\deg(p_{[u]})$, and there
are no decompositions with $\deg(\ell_{T_0})=b_0$, so the right hand side of $(*)$ is $0$.
This proves the second alternative of~\eqref{E:action-us}.
If $2j\geq\deg(p_{[u]})$, then we may rewrite $(*)$ as
\[(p_1\ldots p_u q_{1}\ldots q_{v})\cdot\TO_\beta=\sumsub{T_0\sqcup\cdots\sqcup T_k=[v]\\ \deg(p_{[u]})+\deg(q_{T_0})=b_0,\,\deg(q_{T_i})=b_i} 
\bigl(p_{[u]}q_{T_0}\bigr)
\bigl(q_{T_1}\cdot\TO_{(0,b_1)}\bigr)\ldots\bigl(q_{T_k}\cdot\TO_{(0,b_k)}\bigr)\,.\]
By~\eqref{E:action-GR}, the right hand side of this equation can be written 
in the form
\[\sum_{s\in S_v} c_s  p_1\cdots p_u q_{s(1)}\cdots q_{s(v)}\]
for some choice of scalars $c_s$. This proves the second alternative of~\eqref{E:action-us}. 

We now prove the converse. Fix $j$ and an element $\phi\in\field S_n$ satisfying~\eqref{E:action-us}. This implies that $\phi\in\pppn$, by Theorem~\ref{T:charpeak}. To show that $\phi\in\pppn^j$ we may verify that
 $\pi^{j+1}(\phi)=0$ (Definition~\ref{D:peakideals}).

Choose $j+1$ Lie polynomials $p_0,\ldots,p_{j}$ of degree $2$ and arbitrary homogeneous Lie polynomials $q_1,\ldots,q_v$. Then $2j<\sum_{i=0}^j\deg(p_i)$, so by the second alternative of~\eqref{E:action-us},
\[(p_0\cdots p_j q_1\cdots q_v)\cdot\phi =0\,.\]
On the other hand, applying Lemma~\ref{L:action-pi} repeatedly we obtain
\[(p_0\cdots p_j q_1\cdots q_v)\cdot\phi=(p_0\ldots p_j)\bigl((q_1\ldots q_v)\cdot\pi^{j+1}(\phi)\bigr)\,.\]
Therefore,
\[(p_0\ldots p_j)\bigl((q_1\ldots q_v)\cdot\pi^{j+1}(\phi)\bigr)=0\,,\]
and hence
\[(q_1\ldots q_v)\cdot\pi^{j+1}(\phi)=0\,,\]
since $p_0\ldots p_j$ is a non-zero element of $T(V)$. Now, $q_1\ldots q_v$ is an arbitrary Lie monomial and these span $T(V)$. Since $\dim V\geq n$, the action of $\field S_n$ on $V^{\otimes n}$ is faithful. We conclude that
$\pi^{j+1}(\phi)=0$, which completes the proof.
\end{proof}

Note that if the second alternative of~\eqref{E:action-us} is satisfied then
$2j+2\leq \sum_{i=1}^u \deg(p_i)$, because each $\deg(p_i)$ is even.
In particular, if $j=\ipartn$, the second alternative of~\eqref{E:action-us} is satisfied only when $n<\sum_{i=1}^u \deg(p_i)$, so Theorem~\ref{T:charpeakideal} reduces in this case to Theorem~\ref{T:charpeak}.

On the other hand, the case $j=0$ of Theorem~\ref{T:charpeakideal} recovers Schocker's result~\eqref{E:action-Schocker}. If $\phi\in\pppn^0$,
the theorem implies that $(\ell_1\ldots \ell_k)\cdot\phi=0$ whenever $\ell_1$ is an even Lie polynomial, by the second alternative of~\eqref{E:action-us}.
 Conversely,  suppose an element $\phi\in\SolA$ is such that $(\ell_1\ldots \ell_k)\cdot\phi=0$ whenever $\ell_1$ is an even Lie polynomial. Consider now Lie polynomials $p_1,\ldots,p_u$, $q_1,\ldots,q_v$ as in Theorem~\ref{T:charpeakideal}. If $u=0$, then the first alternative of~\eqref{E:action-us} is satisfied, by~\eqref{E:action-GR}. If $u\geq 1$ then the second alternative is satisfied (with $j=0$), by hypothesis.
The theorem then says that $\phi\in\pppn^0$.

\begin{example}\label{Ex:charpeakideal}
Let $a,b,c,d\in V$ and consider the Lie polynomials $p_1=[a,b]$, $p_2=[c,d]$, and $p_3=a$. We have
\[p_1p_2p_3=abcda-abdca-bacda+badca\,.\]
Consider the element
\[P_{\{4\}}=12354+12453+13452+23451\in\ppp{5}\,.\]
By~\eqref{E:defpi}, $\pi(P_{\{4\}})=P_{\{2\}}\in\ppp{3}$ and $\pi^2(P_{\{4\}})=0$, so 
$P_{\{4\}}\in\ppp{5}^1$.
The action of $P_{\{4\}}$ on $p_1p_2p_3$ may be explicitly calculated as follows:
\begin{align*}
(p_1p_2p_3)\cdot 12354&= abcad-abdac-bacad+badac\\
(p_1p_2p_3)\cdot 12453&= abdac-abcad-badac+bacad\\
(p_1p_2p_3)\cdot 13452&= acdab-adcab-bcdaa+bdcaa\\
(p_1p_2p_3)\cdot 23451&= bcdaa-bdcaa-acdab+adcab\,.
\end{align*}
It follows that
\[(p_1p_2p_3)\cdot P_{\{4\}}= 0\,,\]
in agreement with Theorem~\ref{T:charpeakideal}.
\end{example}



\begin{thebibliography}{10}

\bibitem{ABN}
M. Aguiar, N. Bergeron, and K. Nyman, 
\emph{The peak algebra and the descent algebras of types B and D}, 
to appear in Transactions of the American Mathematical Society.

\bibitem{BB92a}
F. Bergeron and N. Bergeron, \emph{A decomposition of the
descent algebra of the hyperoctahedral group, I}, J. Algebra 148
(1992), no.~1, 86--97.


 \bibitem{BB92b}
F. Bergeron and N. Bergeron,
\emph{Orthogonal idempotents in the descent algebra of $B\sb n$ and applications},
 J. Pure Appl. Algebra  79  (1992), no. 2, 109--129.

\bibitem{B92}
N. Bergeron, \emph{A decomposition of the
descent algebra of the hyperoctahedral group, II}, J. Algebra 148
(1992), no.~2, 98--122.

\bibitem{BHT}
N. Bergeron, F. Hivert, and J.-Y. Thibon,
\emph{The peak algebra and the Hecke-Clifford algebras at $q=0$}, 2004.
  {\tt  math.CO/0304191}
   
 \bibitem{BMSW}
 {\sc N.~Bergeron, S.~Mykytiuk, F.~Sottile, and S.~van Willigenburg}, {\em
   Shifted quasi-symmetric functions and the Hopf algebra of peak
 functions},
  Discrete Math. 246 (2002) 57--66.

\bibitem{GReu}
A. M. Garsia and C. Reutenauer, \emph{A decomposition of Solomon's descent algebra},
 Adv. Math. 77 (1989), no.~2, 189-262.

\bibitem{Hsi}
S. K. Hsiao, \emph{Structure of the peak Hopf algebra of quasi-symmetric
functions}, 2002.

\bibitem{KLT} D. Krob, B. Leclerc and J.-Y. Thibon,
\emph{Noncommutative symmetric functions. II. Transformations of alphabets},
 Internat. J. Algebra Comput.  7  (1997), no. 2, 181--264.

\bibitem{Lam} T. Y. Lam,
\emph{A first course in non-commutative rings},
Graduate Texts in Mathematics 131, Springer-Verlag,1991.

\bibitem{Lod}
Jean-Louis Loday, \emph{Cyclic homology},
Grundlehren der Mathematischen Wissenschaften, 301. 
Springer-Verlag, Berlin, 1998. xx+513 pp. 

\bibitem{Mah01}
S. Mahajan, \emph{Shuffles on Coxeter groups}, 2001. {\tt math.CO/0108094}

\bibitem{MalReu}
C. Malvenuto and C. Reutenauer, \emph{Duality between
  quasi-symmetric functions and the {S}olomon descent algebra}, J. Algebra
  177 (1995), no.~3, 967--982.


\bibitem{Nym}
K. Nyman, \emph{The peak algebra of the symmetric group},
 J. Algebraic Combin.  17 (2003), 309--322.

\bibitem{Reu}
Christophe Reutenauer, \emph{Free {L}ie algebras}, The Clarendon Press Oxford
  University Press, New York, 1993, Oxford Science Publications. 
  
\bibitem{Sch} M. Schocker, 
\emph{The peak algebra of the symmetric group revisited}, to appear in Adv. in Math.  {\tt math.RA/0209376}
 
 \bibitem{Sch2} M. Schocker, 
\emph{The descent algebra of the symmetric group}, to appear in Proc. of ICRA X, Fields Institute, Toronto.

\bibitem{SlSeek}
N. J.~A.~Sloane, \emph{An on-line version of the encyclopedia of integer
  sequences}, Electron. J. Combin. 1 (1994), Feature 1, approx. 5 pp. (electronic);
{\tt http://akpublic.research.att.com/$\sim${}njas/sequences/ol.html}.

\bibitem{Sol}
Louis Solomon, \emph{A Mackey formula in the group ring of a Coxeter group},
 J. Algebra  41  (1976), no. 2, 255--264.

\bibitem{St86}
R.~P. Stanley, \emph{Enumerative combinatorics. {V}ol. 1}, Cambridge
  University Press, Cambridge, 1997, With a foreword by Gian-Carlo Rota.

\bibitem{Ste}
John~R. Stembridge, \emph{Enriched ${P}$-partitions}, Trans. Amer. Math. Soc.
  349 (1997), no.~2, 763--788.

\end{thebibliography}
\end{document}